\documentclass{amsart}

\usepackage{amsthm, amsfonts, amssymb, amsmath, latexsym, enumerate, times}
\usepackage[latin1]{inputenc}
\usepackage{mathrsfs}
\usepackage{mathtools}
\usepackage{letltxmacro}
\usepackage{enumerate}
\usepackage{array}

\usepackage{comment}
\usepackage{todonotes}

\LetLtxMacro\orgvdots\vdots
\LetLtxMacro\orgddots\ddots

\newcounter{rowcntr}[table]
\renewcommand{\therowcntr}{\roman{rowcntr}}

\newcolumntype{N}{>{\refstepcounter{rowcntr}(\therowcntr)}c}

\makeatletter
\DeclareRobustCommand\vdots{%
	\mathpalette\@vdots{}%
}
\newcommand*{\@vdots}[2]{%
	\sbox0{$#1\cdotp\cdotp\cdotp\m@th$}%
	\sbox2{$#1.\m@th$}%
	\vbox{%
		\dimen@=\wd0 %
		\advance\dimen@ -3\ht2 %
		\kern.5\dimen@
		\dimen@=\wd2 %
		\advance\dimen@ -\ht2 %
		\dimen2=\wd0 %
		\advance\dimen2 -\dimen@
		\vbox to \dimen2{%
			\offinterlineskip
			\copy2 \vfill\copy2 \vfill\copy2 %
		}%
	}%
}
\DeclareRobustCommand\ddots{%
	\mathinner{%
		\mathpalette\@ddots{}%
		\mkern\thinmuskip
	}%
}
\newcommand*{\@ddots}[2]{%
	\sbox0{$#1\cdotp\cdotp\cdotp\m@th$}%
	\sbox2{$#1.\m@th$}%
	\vbox{%
		\dimen@=\wd0 %
		\advance\dimen@ -3\ht2 %
		\kern.5\dimen@
		\dimen@=\wd2 %
		\advance\dimen@ -\ht2 %
		\dimen2=\wd0 %
		\advance\dimen2 -\dimen@
		\vbox to \dimen2{%
			\offinterlineskip
			\hbox{$#1\mathpunct{.}\m@th$}%
			\vfill
			\hbox{$#1\mathpunct{\kern\wd2}\mathpunct{.}\m@th$}%
			\vfill
			\hbox{$#1\mathpunct{\kern\wd2}\mathpunct{\kern\wd2}\mathpunct{.}\m@th$}%
		}%
	}%
}
\makeatother

\usepackage{array}
\usepackage{comment}
\usepackage{paralist}
\usepackage{xcolor}

\newtheorem{introtheorem}{Theorem}
\newtheorem{theorem}{Theorem}[section]
\newtheorem{lemma}[theorem]{Lemma}
\newtheorem{claim}[theorem]{Claim}
\newtheorem{corollary}[theorem]{Corollary}
\newtheorem{proposition}[theorem]{Proposition}

\theoremstyle{definition}

\newtheorem{remark}[theorem]{Remark}
\newtheorem{example}[theorem]{Example}

\newcommand{\calL}{{\mathcal L}}

\newcommand{\calO}{{\mathcal O}}

\newcommand{\bbF}{{\mathbb F}}
\newcommand{\bbP}{{\mathbb P}}

\def\geq{\geqslant}
\def\leq{\leqslant}
\def\le{\leqslant}
\def\ge{\geqslant}

\begin{document}
 
\title[Planar Linear Systems Assuming SHGH]{Boundedness Results For Planar Linear Systems Assuming The Segre-Harbourne-Gimigliano-Hirschowitz Conjecture}
 
\author{Ciro Ciliberto}
\address{Dipartimento di Matematica, Universit\`a di Roma Tor Vergata, 
Via O. Raimondo 00173 Roma, Italia}
\email{cilibert@axp.mat.uniroma2.it}

\author{Rick Miranda}
\address{Department of Mathematics, Colorado State University, 
Fort Collins (CO), 80523,USA}
\email{rick.miranda@colostate.edu}
 
\author{Joaquim Ro\'e}
\address{Departament de Matematiques, Universitat Aut\`onoma de Barcelona, 
08193 Bellaterra (Barcelona), Catalunya}
\email{joaquim.roe@uab.cat}

 

\maketitle

\begin{abstract}
	Let $X_n$ be the projective plane blown up at $n \geq 10$ general points.
	In this paper we give several consequences of the Segre-Harbourne-Gimigliano-Hirschowitz Conjecture, that pertain to complete linear systems on $X_n$.
	We begin by classifying such systems $|C|$
	with general irreducible member of genus $g \geq 2$ (up to Cremona equivalence),
	in terms of invariants of the adjoint systems $|C+mK|$.
	We then use this to prove that, for fixed $n \geq 10$ and $g\geq 2$,
	up to the action of the Cremona group,
	there exist finitely many complete linear systems on $X_n$
	whose general member is irreducible of genus $g$.
	Further, there is a function $g\mapsto n(g)$ such that
	every such (effective) system is Cremona equivalent to a system in $X_{n(g)}$.
	The latter result is based on the explicit computation of the
	minimum possible self-intersection of an irreducible linear system
	with given $n$ and $\dim(|C|)$.
	We classify those systems which achieve the minimal self-intersection.
	We also classify the systems with $C^2 \leq 5$, whether or not they have minimal $C^2$ for the given $n$ and dimension.
	We finish by proving several statements concerning systems that are base-point-free,
	and systems that give birational maps to their image.
	
\end{abstract}

\tableofcontents
 
\section*{Introduction} 

Let $X_n$ be the complex projective plane blown up at $n$ general points. 
Let $E_1,\ldots, E_n$ be the exceptional divisors on $X_n$ over the $n$ blown--up points 
and let $H$ be the class on $X_n$ of the pull back of a line of the plane. 
Then a complete linear system $\calL$ on $X_n$ can be written as 
$|aH-\sum_{i=1}^n b_i E_i|$ which we will abbreviate as $|a;b_1,\ldots, b_n|$; 
$a$ will be called the \emph{degree} of the system 
and $b_1,\ldots,b_n$ the \emph{multiplicities},
and the collection of these integers are called the \emph{numerical characters} of the system.
After permuting the blown--up points, 
we may assume that $b_1\geq b_2 \geq \cdots \geq b_n$.  
A linear system $|a;b_1,\ldots, b_n|$ is \emph{Cremona reduced}
(or \emph{Cremona minimal}) if $a\geq b_1+b_2+b_3$; 
it is \emph{$n$-Cremona minimal} if it is Cremona minimal with $b_n \geq 1$.
(This essentially means that the linear system 
is not pulled back from a blowup of the plane at fewer than $n$ points,
and is Cremona minimal.)
Two linear systems are said to be \emph{Cremona equivalent} if there is a Cremona transformation that maps one to the other.

The \emph{self-intersection} of $\calL$ as above is given by $\calL^2 = a^2-\sum_{i=1}^n b_i^2$ and the 
	(arithmetic)
\emph{genus} of $\mathcal L$ is 
$$g(\calL)=\frac {(a-1)(a-2)}2-\sum_{i=1}^n \frac {b_i(b_i-1)}2.$$
Self-intersection and genus are invariant under Cremona transformations.
Cremona minimal linear systems have minimal degree among their Cremona equivalent linear systems, but, if the base points are not general, there are examples of  Cremona minimal and Cremona equivalent linear systems that are not projectively equivalent (see \cite [Thm. 0.2]{CCi}).


We make the following assumptions, while establishing some notation:

\begin{enumerate}[(A)]
\item \label{hypothesis:SHGH} 
The Segre-Harbourne-Gimigliano-Hirschowitz (SHGH) Conjecture
(see \cite{CM}) ;
\item \label{hypothesis:g} 
$\calL$ is a linear system given by $|a;b_1,\ldots,b_n|$ 
satisfying $r:=\dim(\calL)\geq 0$ and $g:=g(\calL)\geq 2$.
\item \label{hypothesis:no-less-than-n-points}
we have $b_1\geq b_2 \geq \cdots \geq b_n\geq 1$ and
for any $(-1)$--curve $E$ on $X_n$, we have $\calL\cdot E \geq 1$,
\end{enumerate}

If $E$ is a $(-1)$--curve on $X_n$ with $\calL\cdot E = 0$,
then after a Cremona transformation we can reduce $n$,
so we call systems satisfying \eqref{hypothesis:no-less-than-n-points}
\emph{$n$-point systems}. 
Thus, a Cremona minimal system is $n$-Cremona minimal
if and only if it is an $n$-point system.

With these assumptions at hand, we can fix further notations and definitions, 
and make some initial observations:
\begin{enumerate}[(A)]
\setcounter{enumi}{3}
\item \label{hypothesis:Kn}
We let $K_n$ denote the canonical class on $X_n$, 
whose linear system is $|{-3;(-1)^n}|$
	(we use exponential notation for repeated multiplicities)
\item \label{hypothesis:negcurves} 
If $E$ is an irreducible curve on $X_n$ with $E^2<0$,
then $E$ is a $(-1)$--curve (a smooth rational curve with $E^2=-1$).
\item \label{hypothesis:cohomology}
If $F$ is an effective divisor, and is either nef or reduced, then
we have $h^i(X_n,F) = 0$ for $i=1,2$, 
and $\dim(|F|) = F^2-g(F)+1$.
\end{enumerate}
We will see below (Propositions \ref {prop:g2} and \ref{prop:adj})
that the following two statements (G) and (H)
follow from (A), (B), and (C) also:
\begin{enumerate}[(A)]
\setcounter{enumi}{6}
\item \label{hypothesis:nobase} 
The general curve $C\in \calL$ is irreducible.
\item \label{hypothesis:adjoint_nofix}
The adjoint system $|C+K_n|$, of dimension $g(\calL)-1$, has no fixed components;
its general member will be denoted by $C'$.
\end{enumerate}

The SHGH Conjecture (A) is a strong assumption: 
we may view this work as providing some consequences of that conjecture.
It is well known that it implies
\eqref{hypothesis:negcurves} and \eqref{hypothesis:cohomology}, 
where the dimension formula 
follows from cohomology vanishing by the Riemann-Roch theorem;
see \cite{CM,CHMR}.
The vanishing statement for reduced effective divisors 
is known as \emph{Segre's Conjecture.}
It also implies that every linear system $\calL$ 
satisfying \eqref{hypothesis:g} and  \eqref{hypothesis:no-less-than-n-points} is nef.

	
Let us comment further on the assumption \eqref{hypothesis:g}. 
According to the SHGH Conjecture and its consequences 
\eqref{hypothesis:negcurves} and \eqref{hypothesis:cohomology}, 
if $C$ is an irreducible rational curve 
then either it is a $(-1)$--curve or $C^2\geq 0$ 
and $|C|$ is base point free of dimension $C^2+1$.
The Cremona minimal systems Cremona equivalent to such a $|C|$ 
are well known by \cite [Chapt. 5]{Ca}. 
If $C$ is irreducible of genus 1, 
then $|C|$ has dimension $C^2\geq 0$ 
and the Cremona minimal systems Cremona equivalent to $|C|$ are again well known 
(see  \cite [Chapt. 5]{Ca} and \cite [Prop. 10.6] {CCi});
they are $|3;1^n|$ and $|4;2^2|$.
So, the cases $g\leq 1$ being well known, we are assuming $g\geq 2$. 

Our first finiteness result is the following.
\begin{introtheorem}\label{thm:main} 
Assume the SHGH conjecture holds.
Given an integer $g\geq 2$,
there are finitely many 
sets of numerical characters $|a;b_1,\ldots,b_n|$ of 
effective Cremona minimal
linear systems $\calL$
for which $g(\calL)=g$.
\end{introtheorem}

We pose the question of whether the finiteness noted above holds
even for systems where the points are not in general position,
or even for $g=1$, of dimension at least two.
(The Halphen pencils give counterexamples of degree a multiple of $3$ for $g=1$.)
We are not aware of any examples that would counter this.

Next we consider how the questions at hand depend on
the number $n$ of blown-up points. 
Taking into account that the cases with $n\le 9$,
where the anticanonical divisor is effective, are much better known,
in the last part of the paper we focus on the cases with $n \geq 10$.
We have the following bound on the self-intersection of the systems in question.

\begin{introtheorem}\label{thm:bound_selfintersection}
Assume that the SHGH Conjecture holds,
and that $n \geq 10$.
Let $|C|$ be an effective complete $n$-point linear system on $X_n$
of genus $g \geq 2$ and dimension $r \geq 0$.
(These essentially comprise assumptions 
\eqref{hypothesis:SHGH}, \eqref{hypothesis:g}, 
and \eqref{hypothesis:no-less-than-n-points} above.)
Then:
\begin{itemize}
\item[(a)] If the adjoint system $|C+K_n|$ is a pencil or is composed with a pencil,
then $|C|$ is Cremona equivalent to either
\begin{itemize}
\item[(i)] $|6;2^8,1^{n-8}|$ with $8 \leq n \leq 11$, in which case $g=2$ and $C^2=12-n$;
\item[(ii)] $|g+2;g,1^{n-1}|$ with $1 \leq n \leq 3g+6$, 
in which case $C^2 = 4g-n+5=(n+4r)/3 - 2$
\end{itemize}
\item[(b)] If $g \geq 3$, $n+r=2h$ is even, 
and the adjoint system $|C+K_n|$ has irreducible general member,
then $C^2 \geq h+r-5$.
\item[(c)] If $g \geq 3$, $n+r=2h+1$ is odd, 
and the adjoint system $|C+K_n|$ has irreducible general member,
then $C^2 \geq \lceil (6h-7)/5 \rceil +r-5$.
\end{itemize}
\end{introtheorem}

Note that if $g=2$ then the adjoint system $|C+K_n|$ is a pencil ($|C|$ is hyperelliptic), so, by observation \eqref{hypothesis:cohomology}, cases (a)--(c) cover all possibilities
in Theorem \ref{thm:bound_selfintersection}.

The proof of Theorem \ref{thm:bound_selfintersection}
occupies Section \ref{sec:bound_selfintersection}
with Theorem \ref{thm:hyperelliptic},
Theorem \ref{thm:irreduciblen+reven}, and
Theorem \ref{thm:irreduciblen+rodd}
covering the cases (a), (b), and (c) respectively.
The reader will notice that in the hyperelliptic case (a),
the self-intersection is determined by $n$ and $r$,
while in cases (b) and (c) we only have a lower bound.
Those lower bounds are sharp, and in Theorems
\ref{thm:irreduciblen+reven} and \ref{thm:irreduciblen+rodd}
we also present the systems that achieve the minimum,
in all cases except the lowest possible values for $n$ and $r$.

In Proposition \ref{prop:minimalatmostfive} and Theorem \ref{C2leq5}
we determine all systems with $C^2 \leq 5$,
assuming the hypotheses of (A), (B), and (C) above,
and that $n \geq 10$, $r\geq 2$
(although we present some results for $r=0,1$ as well).

%
%
%

The paper is organized as follows. 
In section \ref{sec:cr} we classify Cremona--minimal linear systems, 
following \cite{CCi}, which is then used in section \ref{sec:finiteness} 
to prove Theorem \ref{thm:main}. 
Section \ref{sec:bound_selfintersection} is devoted to 
the study of linear systems of minimal selfintersection, 
and culminates in the proof of Theorem \ref{thm:bound_selfintersection}.
More precisely, we compute the exact minimum possible value of $C^2$ 
for a given $n$ and $r$, which surprisingly (as noted above) 
behaves differently depending on the parity of $n+r$.
Finally, in section \ref{sec:final} we collect a few consequences of our results,
which include the classification of $C^2 \leq 5$ systems noted above,
as well as several remarks including the gaps in the possible self-intersections
(given $(n,r)$),
the birationality of the maps given by systems with minimal self-intersection,
and several interesting examples with the minimum $n=10$ cases.

\section{Basic consequences of the SHGH conjecture}

As said above, we are going to assume that the SHGH conjecture holds. 
In this section we point out a number of basic consequences of this assumption.

\begin{lemma}\label{lem:fix} 
Assume that the SHGH conjecture holds. 
Let $\calL$ be a linear system on $X_n$ of dimension $r\geq 0$ 
and let $F$ be an irreducible curve that is a fixed component for $\calL$. 
Then
$$
F^2=g(F)-1.
$$
\end{lemma}

\begin{proof} 
Since $F$ is a fixed component for $\calL$, 
one has $\dim(|F|)=0$,
and since it is reduced,
we have (using consequence (\ref{hypothesis:cohomology}) of SHGH)
that its higher cohomology vanishes. Therefore
$$
0=\dim(|F|)=\frac {F\cdot (F-K_n)}2=F^2-g(F)+1
$$
and the assertion holds. 
\end{proof}

\begin{lemma}\label{lem:fix} 
Assume that the SHGH conjecture holds. 
Let $\calL$ be a linear system on $X_n$  of dimension $r\geq 0$ 
and let $F$ be an irreducible curve that is a fixed component for $\calL$. 
Then if $D\neq F$ is an irreducible curve such that 
$D$ is contained in a member of $\calL$, then $D\cdot F=0$.
\end{lemma}

\begin{proof} 
By the assumption, we must have $\dim(|D+F|)=\dim(|D|)$.
Since the SHGH conjecture holds,
and both $D$ and $D+F$ are reduced,
we have that both their higher cohomology vanishes
(using (\ref{hypothesis:cohomology})),
and so
\begin{align*}
\frac {D\cdot (D-K_n)}2&=\dim(|D|)=\dim(|D+F|)=\frac {(D+F)\cdot (D+F-K_n)}2=\\
&=\frac {D\cdot (D-K_n)}2+\frac {F\cdot (F-K_n)}2+F\cdot D
\end{align*}
and the assertion follows. 
\end{proof}

\begin{proposition}\label{prop:-1} 
Assume that the SHGH conjecture holds. 
Let $\calL$ be a linear system on $X_n$  of dimension $r\geq 0$ 
and let $F$ be a $(-1)$--curve that is a fixed component for $\calL$.  
Then $\calL=aF+\calL'$, with $a$ a positive number, 
$F$ not a fixed component of $\mathcal L'$ 
and $\calL'\cdot F=0$. 
In particular $\calL$ is not nef. 

Conversely, if $\calL$ is not nef, then there is a $(-1)$--curve in the fixed part of $\calL$. 
\end{proposition}

\begin{proof} 
The first part of the assertion is an immediate consequence of Lemma \ref {lem:fix} 
and if $\calL=aF+\calL'$ with $F$ a $(-1)$--curve, 
then $\calL\cdot F=-a<0$, so $\calL$ is not nef. 
Conversely, if $\calL$ is not nef, 
there is some irreducible curve $F$ with $F^2<0$ 
such that $\calL\cdot F<0$ and $F$ is fixed for $\calL$. 
Since SHGH holds, consequence \eqref{hypothesis:negcurves} implies that 
$F$ must be a $(-1)$--curve and the final statement follows. 
\end{proof}

We will now focus on the case where
$\calL$ is a nef a linear system on $X_n$ of dimension $r\geq 0$, 
so that there is no $(-1)$--curve in the fixed part of $\calL$. 

\begin{lemma}\label{lem:g1} 
Assume that the SHGH conjecture holds. 
Let $\calL$ be a nef linear system on $X_n$ of dimension $r\geq 0$ 
and let $F$ be an irreducible curve that is a fixed component for $\calL$. 
If there is an integer $a\geq 2$ such that $aF$ is fixed for $\calL$ 
then $F^2=0$, $g(F)=1$, 
and $|F|$ is Cremona equivalent to $|3;1^9, 0^{n-9}|$.
\end{lemma}

\begin{proof}
Since $\calL$ is nef, by Proposition \ref{prop:-1} 
the curve $F$ is not a $(-1)$-curve,
and hence $F^2\geq 0$ by consequence (\ref{hypothesis:negcurves}) of SHGH.
Using Lemma \ref{lem:fix}, we have $F^2 = g(F)-1$;
then we have 
$$
0=\dim(|2F|)= F(2F-K)=3F^2-F(F+K)=3F^2-2(g(F)-1)=F^2,
$$
so that $g(F)=1$.
We conclude that $|F|$ is Cremona equivalent to $|3;1^9, 0^{n-9}|$ 
(see  \cite [Chapt. 5]{Ca}). 
\end{proof}

\begin{proposition}\label{prop:g11} 
Assume that the SHGH conjecture holds. 
Let $\calL$ be a nef linear system on $X_n$ of dimension $r\geq 0$ 
and let $F$ be an irreducible curve that is a fixed component for $\calL$. 
If $F^2=0$ (and then $g(F)=1$) 
there is a positive integer $m$ such that $\calL=\{mF\}$, or, equivalently, 
$\calL$ is Cremona equivalent to $|3m;m^9, 0^{n-9}|$ and $g(\calL)=1$.
\end{proposition}

\begin{proof} 
Let $m$ be the maximum integer such that $mF$ is fixed for $\calL$. 
By Lemma \ref {lem:fix}, we have $(\mathcal L-mF)\cdot F=0$, 
and so $\mathcal L\cdot F=0$ since we are assuming $F^2=0$.
Suppose $\calL\neq \{mF\}$;
%
then there is some irreducible curve $D\neq F$ 
such that $mF+D$ is contained in a member of $\calL$. 
Then $D$ is not a $(-1)$--curve because $\calL$ is nef 
(by Proposition \ref {prop:-1}) 
and $D\cdot F=0$ by Lemma \ref{lem:fix}. 
Therefore $D^2\geq 0$, 
and if $D^2>0$ we have a contradiction to the index theorem. 
If $D^2=0$, then $D$ and $F$ would be $\mathbb Q$--numerically equivalent, 
which is also impossible. 
The assertion follows.
\end{proof}

\begin{proposition}\label{prop:g2} 
Assume that the SHGH conjecture holds.
Let $\calL$ be a nef linear system on $X_n$ of dimension $r\geq 0$ 
that is not Cremona equivalent to $|3m;m^9, 0^{n-9}|$ 
for some positive integer $m$. 
Then either $\calL$ has no fixed component 
or $r=0$ and the unique curve in $\calL$ is irreducible.
\end{proposition}

\begin{proof}
Suppose first $r>0$ and that there is a fixed irreducible curve $F$ for $\calL$.
By Lemma \ref {prop:-1} and Proposition \ref{prop:g11}, one has $F^2>0$.
Let $|M|$ be the movable part of $\calL$, so that $M^2\geq 0$. 
By Lemma \ref {lem:fix}, one has $M\cdot F=0$, 
which contradicts the index theorem. 

Hence we must have $r=0$. 
Notice that no multiple of $F$ can be fixed in $\calL$ by Lemma \ref {lem:g1}. Suppose $\calL\neq \{F\}$. 
Then there is a curve $D\neq F$ such that 
$F+D$ is contained in the unique curve of $\calL$ and $D\cdot F=0$. 
By the same reasons as above we must have $D^2>0$, 
which again contradicts the index theorem.
\end{proof}

\begin{proposition}\label{prop:fin} 
Assume that the SHGH conjecture holds. 
Let $\calL$ be a nef linear system of dimension $r\geq 1$ on $X_n$. 
Then either the general curve in $\calL$ is irreducible, 
or $\calL$ is composed with a pencil $|P|$ of rational curves, with $P^2=0$. 
In this case $\calL$ is Cremona equivalent to $|d;d,0^{n-1}|$ 
for some integer $d\geq 2$ 
and $g(\calL)=1-d$. 
\end{proposition}

\begin{proof} 
By Proposition \ref {prop:g2}, 
we know that $\calL$ has no fixed component. 
If the general curve of $\calL$ is not irreducible, 
then there is a pencil $|P|$ such that $\calL=|dP|$, with $d\geq 2$. 
First of all we have
$$
1=\dim (|P|)=\frac {P\cdot (P-K_n)}{2} = P^2-g(P)+1
$$
that yields $P^2=g(P)$. 
On the other hand we have $\dim(|2P|)=2$, which implies 
$$
2=\dim(|2P|)=P(2P-K_n)=2P^2-P\cdot K_n=3P^2-2g(P)+2
$$
so that
$$
3g(P)=3P^2=2g(P)
$$
which gives $g(P)=0=P^2$. 
The final assertion follows from classical results (see \cite [Chapt. 5]{Ca}).  \end{proof}

Propositions \ref{prop:g2} and \ref{prop:fin}
show that hypotheses (\ref{hypothesis:SHGH}), (\ref{hypothesis:g}), and (\ref{hypothesis:no-less-than-n-points})
imply (\ref{hypothesis:nobase}).
We first note that the three hypotheses imply that $\calL$ is nef.
Then Proposition \ref{prop:g2} gives that either $r=0$ and the single member of $\calL$ is irreducible, or $r\geq 1$ and $\calL$ has no fixed components.
Then Proposition \ref{prop:fin} finishes the argument, since $g(\calL) \geq 2$.

The following completes the argument for (\ref{hypothesis:adjoint_nofix}).

\begin{proposition}\label{prop:adj} 
Assume that the SHGH conjecture holds. 
Let $|C|=\calL$ be a nef linear system on $X_n$ 
satisfying hypothesis \eqref {hypothesis:g},
i.e., $r \geq 0$ and $g(\calL) \geq 2$. 
Then the adjoint system $|C + K_n|$, 
of dimension  $g(\calL)-1$, 
has no fixed components.
\end{proposition}

\begin{proof}
First note that, since $g(\calL) \geq 2$,
$\calL \neq |3m;m^9,0^{n-9}|$,
so that by Proposition \ref {prop:g2},
either $\calL$ has no fixed component, 
or $\calL$ consists of a single irreducible curve.
However if $\calL$ has no fixed component, then $r \geq 1$,
so that Proposition \ref{prop:fin} applies,
and we conclude that either $\calL$ is composed with a rational pencil
(which is forbidden by $g(\calL) \geq 2$)
or the general curve in $\calL$ is irreducible.
We conclude that in any case $\calL$ has an irreducible curve $C$ as a member.
By looking at the exact sequence
$$
0\longrightarrow \mathcal O_{X_n}(K_n)
\longrightarrow \mathcal O_{X_n}(C+K_n)
\longrightarrow \omega_C\longrightarrow 0
$$
and taking into account that 
$h^i(X_n, \mathcal O_{X_n}(K_n))=0$ for $0\leq i\leq 1$ 
and $h^0(C,\omega_C)=g(\calL)$,
we see that the adjoint system $|C + K_n|$ has dimension $g(\calL)-1\geq 1$. Moreover it cuts out the complete canonical system on $C$, 
which has no base points. 
Hence, if there is a fixed irreducible curve $E$ for $|C + K_n|$, 
one has $E\cdot C=0$. 
On the other hand we have 
$0\leq \dim (\calL)=C^2-g(\calL)+1$, 
thus $C^2\geq g(\calL)-1\geq 1$. 
Hence, by the index theorem, one has $E^2<0$, so 
(using consequence \eqref{hypothesis:negcurves})
$E$ is a $(-1)$--curve such that $C\cdot E=0$, 
a contradiction.
\end{proof}

\section{Cremona minimality}\label{sec:cr}

In this section we will classify 
the effective Cremona minimal linear systems $\mathcal L$ of genus $g\ge 2$. 
This has been essentially done in \cite {CCi}, 
but we think it is useful to work out here this particular case. 
We keep all notation and conventions established in the Introduction.

\begin{lemma}\label{lem:a}
Assume that the SHGH Conjecture holds
and let $|C|$ be 
a linear system on $X_n$ of genus $g\geq 2$.
Suppose that the \emph{$m$--adjoint system} $|C+mK_n|$
is effective but not nef.
Then there are $h\leq n$ disjoint $(-1)$--curves $A_1,\ldots, A_h$ on $X_n$ 
such that  $0>(C+mK_n)\cdot A_i=C\cdot A_i-m$ for $1\leq i\leq h$.
Contracting these $(-1)$--curves gives a morphism $f: X_n\longrightarrow S$ 
(where either $S\cong X_{n-h}$ or $S\cong \mathbb F_0$).
If we set  $C_S=f_*(C)$, then
$|C_S+mK_S|$ is effective and nef, and 
$$
h^0(X_n, \mathcal O_{X_n}(C+mK_n))=h^0(S, \mathcal O_{S}(C_S+mK_S)).
$$
\end{lemma}

\begin{proof} 
If $|C+mK_n|$ is effective but not nef, 
there is some irreducible curve $A$ such that $A\cdot (C+mK_n)<0$. 
Since $A^2$ must be strictly negative, 
consequence (\ref{hypothesis:negcurves}) implies that
$A$ must be a $(-1)$--curve. 
Moreover, if $A,B$ are two $(-1)$--curves such that 
$A\cdot (C+mK_n)<0$ and $B\cdot (C+mK_n)<0$, 
one must have $A\cdot B=0$, because otherwise $A+B$  
would be an effective divisor which meets each of its components nonnegatively, therefore nef, contradicting $(A+B)\cdot (C+mK_n)<0$
since $C+mK_n$ is effective. 
The curves $A_i$ in the statement are then all such $(-1)$--curves,
and there are at most $n$ of them since 
they generate a negative definite rank $h$ lattice inside $NS(X_n)$.
This proves the first part of the assertion. 

Contracting the $h$ disjoint $(-1)$--curves $A_1,\ldots, A_h$
gives the morphism $f: X_n\longrightarrow S$ 
where either $S\cong X_{n-h}$ or $S\cong \mathbb F_0$ 
(which can only happen if $h=n-1$). 
Set $e_i:=C\cdot A_i$;
note that $e_i<m$ since $A_i\cdot(C+mK_n) < 0$, and $A_i\cdot K_n = -1$.
Then $C=f^*(C_S)-\sum_{i=1}^h e_iA_i$. 
Moreover $K_n=f^*(K_{S})+\sum_{i=1}^h A_i$. 
So
\begin{equation}\label{eq:a}
C+mK_n\equiv f^*(C_S+mK_{S})+\sum_{i=1}^h (m-e_i)A_i
\end{equation}
which proves the final assertion. 
\end{proof}

\begin{remark}\label{rem:b} 
The preceding Lemma \ref {lem:a} can be rephrased as an explicit description 
of the Zariski decomposition of $C+mK_n$, namely:
$$
C+mK_n\sim P+ A
$$
where 
\begin{itemize}
\item $P=f^*(C_S+mK_{S})$ is nef, 
\item $A=\sum_{i=1}^h (m-e_i)A_i $ is effective 
with negative definite intersection matrix,
\item $P\cdot A_i=0$ for every irreducible component $A_i$ of $A$, and
\item  $h^0(X_n, \mathcal O_{X_n}(C+mK_n))=h^0(X_n, \mathcal O_{X_n}(P))$.
\end{itemize}
\end{remark}

Next, given a linear system $\calL=|C|$, 
we introduce the following invariants:
\begin{itemize}
\item $m:=m(\calL)$ is the minimum (positive) integer such that 
$|C+mK_n|\neq \emptyset$ and $|C+(m+1)K_n|= \emptyset$;
\item  $\alpha:=\alpha(\calL)=\dim (|C+mK_n|)$.
\end{itemize} 
With these invariants in hand we will be able to 
provide the desired classification of Cremona minimal linear systems.

\begin{theorem}\label{thm:mm} 
Assume that the SHGH conjecture holds.
Let $\calL$ be an effective $n$-point linear system 
with $g(\calL) \geq 2$,
and let $m=m(\calL)$, $\alpha=\alpha(\calL)$ be as defined above.
Then $\calL$ is Cremona equivalent to 
one of the following $n$-Cremona minimal systems:
\begin{enumerate}[(i)]
\item \label{lastadjoint:plurianticanonical}$|3m;m_1,\ldots,m_{n}|$ with $m\geq m_1\geq \cdots \geq m_n\geq 1$;\\
this case can occur for any $m\geq 1$, but only occurs if $\alpha=0$.
\item \label{lastadjoint:plurianticanonical_F0}$|3m+e;m+ e,m+e,m_3,\ldots, m_n|$ with $m-e\geq m_3\geq \cdots \geq m_{n}\geq 1$, $m>e>0$;\\
this case occurs for any $m \geq 2$ but only if $\alpha= 0$.
\item \label{lastadjoint:pencilP}$|3m+\alpha;m+ \alpha,m_2,\ldots, m_n|$ with $m\geq m_2\geq \cdots \geq m_{n}\geq 1$;\\
this case occurs for any $m \geq 1$, and any $\alpha\geq 1$.
\item \label{lastadjoint:F_0} $|3m+\alpha+e;m+\alpha+e, m+e, m_3, \ldots, m_{n}|$ with  $m-e\geq m_3\geq \cdots \geq m_{n-1}\geq 1$, $m>e>0$;\\
this case occurs for any $m \geq 2$, and any $\alpha\geq 1$.
\item \label{lastadjoint:veronese} $| 3 m + \left\lfloor \frac {\alpha} 2\right\rfloor ; m_1 , \ldots, m_{n} |$ , with $m\geq m_1\geq \cdots \geq m_{n}\geq 1$;\\
this case occurs for any $m \geq 1$ but only if $\alpha\in\{2,5\}$.
\item \label{lastadjoint:big_even} 
$|3m+\frac \alpha 2;m-1+\frac\alpha2,m_2,\ldots , m_n|$, 
with $ m\geq m_2\cdots \geq m_n\geq 1$;\\
this case occurs for any $m \geq 1$ but only if $\alpha\geq 4$ is even.
\item \label{lastadjoint:big_odd}
$|3m +\frac {\alpha+1}2+e;m+\frac {\alpha-1}2+e,m+e+1, m_3, \ldots ,m_n|$ with $m-e\geq m_3\geq \cdots \geq m_n\geq 1$ and $m>e\ge 0$;\\
this case occurs for any $m \geq 1$ but only if $\alpha\geq 3$ is odd.
\end{enumerate}
\end{theorem} 

\begin{proof} 
The hypotheses are essentially
\eqref{hypothesis:SHGH}, \eqref{hypothesis:g}, 
and \eqref{hypothesis:no-less-than-n-points} of the Introduction.
The proof  follows from the analysis of three different cases, according to the behavior of the nef part of $|C +mK_n|$. 
Indeed, the Zariski decomposition above of $C +mK_n$ as $P+A$, 
where $P$ is nef and 
$A$ is the sum of multiples of a certain number $h$ of disjoint $(-1)$--curves,
gives rise to three possible cases:\\
\begin{inparaenum}[(a)]
\item \label{case:fixed} $ P= 0$;\\
\item \label{case:pencil} $ P >0$ and $P^2 = 0$;\\
\item \label{case:general} $ P> 0$ and $P^2 >0$.
\end{inparaenum}

Suppose first we are in case \eqref{case:fixed},
which implies that $\alpha = 0$.
Let us blow down the curves in 
$A=\sum_{i=1}^h e_i A_i$
(ordered so that the $e_i$ are non-increasing)
via a morphism $f: X_n\longrightarrow S$,
where either $S\cong X_{n-h}$ or $S\cong \mathbb F_0$;
then $|f_*(C)| =  |-mK_S|$. 
If $S\cong X_{n-h}$ then 
$|f_*(C)| = |3m;m^{n-h}|$ so that $|C|$ is Cremona-equivalent to 
$|3m;m^{n-h},e_1,\ldots,e_h|$
which gives \eqref{lastadjoint:plurianticanonical}.
Note that $n-h$ is the number of $m_i$ equal to $m$, and
we know that $m_n=e_h\geq 1$ because of 
assumption \eqref{hypothesis:no-less-than-n-points}.

If $S\cong \mathbb F_0$  (so $h=n-1$), 
let $F,G$ be the classes of the two rulings of $\mathbb F_0$. 
Then $|f_*(C)|$ is Cremona equivalent 
to a linear system of curves in $|2mF+2mG|$ on $\mathbb F_0$ 
with points of multiplicities 
$e_1,\ldots, e_h$ with $m>e_i>0$ (as in the previous lemma).
Set $e=m-e_1$.
We note that the one-point blowup of $\bbF_0$ is isomorphic to $X_2$,
and the linear system $|aF+bG-cA_1|$ on $\bbF_0$
corresponds to the system $|a+b-c;a-c,b-c|$ on $X_2$.
Hence the system $|2mF+2mG-e_1A_1|$ corresponds to 
$|4m-e_1;2m-e_1,2m-e_1|=|3m+e;m+e,m+e|$ on $X_2$.
(Interpreting $\bbF_0$ as a smooth quadric in $\bbP^3$, 
this is achieved geometrically by projecting down to $\bbP^2$ 
from the point $p$ of largest multiplicity $e_1$.
This contracts the two rulings through $p$ to the two points in $X_2$.)
Since $n=h+1$ we end up with a system of type \eqref{lastadjoint:plurianticanonical_F0} with $m_i=e_{i-1}$ for $i\geq 3$.
	
Suppose next we are in case \eqref{case:pencil}. 
Since $|P+K_n|=\emptyset$, 
we see that $|P|$ is composed 
with curves of an irreducible base point free pencil $|L|$ of rational curves, 
i.e., $P\sim \alpha L$. 
We blow down $A$ and then, successively, 
all other $(-1)$--curves $B$ such that $P\cdot B=0$.
The latter blowing down is not uniquely determined, 
but it eventually provides a morphism $f: X_n\longrightarrow S$, 
where $g: S\longrightarrow \bbP^1$ is a relatively minimal fibration $|F|$
with smooth rational fibres, 
and $|L|$ is the pull back to $X_n$ of $|F|$. 
Then $S\cong \bbF_k$ where either $k=1$ or $k=0$, 
because of assumption \eqref{hypothesis:negcurves}. 
In the former case we denote by $G$ the $(-1)$--curve on $\bbF_1$, 
in the latter we denote by $|G|$ the ruling of $\bbF_0$ different from $|F|$. 
We have
$$
D:= f_*(C) =-mK_S +\alpha F =2mG +(m(2+k) +\alpha) F, 
$$
(since on $\bbF_1$ we have $K=2G+3F$
and on $\bbF_0$ we have $K=2G+2F$)
and $f_*(\calL)$ is a sublinear system of $|D|$ of curves 
with $h$ points of multiplicity $e_i\leq m$. 
If $k=1$ we 
blow down $G$ to a point of multiplicity $D\cdot G=m+\alpha$ to
end up in case \eqref{lastadjoint:pencilP}, 
whereas if $k=0$ we again project from
the point of largest multiplicity $e_1$ to obtain two points
of multiplicity $m+e$ and $m+\alpha+e$ where $e=m-e_1$, so we
are in case \eqref{lastadjoint:F_0}.
Note that in the case $e=0$ the non-uniqueness of $B$ manifests itself: 
via a flip the surface $S=\bbF_0$ can be replaced by $S=\bbF_1$ 
and we fall back to case \eqref{lastadjoint:pencilP} with $m_2=m$.
Thus we may assume $e>0$.

Finally, we deal with case \eqref{case:general}. Since $|P+K_n|=\emptyset$, we see that $|P|$ is a base point free linear system of dimension
$\alpha = P^2 + 1 >1$ of rational curves. Then there is a morphism $f: X_n\longrightarrow S$, with either $S\cong \mathbb P^2$, or $S\cong \mathbb F_1$, or $S\cong \mathbb F_0$ such that $|P'|=f_*(|P|)$ is one of the following:\\
\begin{inparaenum}[(1)]
\item $S\cong \mathbb P^2$ and $|P'|=|\mathcal O_{\mathbb P^2}(j)|$, with $1\leq j\leq 2$, so that $\alpha=2$ in the former case, $\alpha=5$ in the latter;
\item $S\cong \mathbb F_1$ and $|P'|$ is of the form $|d;d-1|$, with $d\geq 2$, in which case $\alpha=2d\geq 4$ is even;\\
\item $S\cong \mathbb F_0$ and $|P'|$ is of the form $|F+dG|$, with $d\geq 1$, where $|F|, |G|$ are the two rulings of $\mathbb F_0$,
in which case $\alpha=2d+1\geq 3$ is odd.
\end{inparaenum}
These three cases clearly lead to cases (v), (vi) and (vii) respectively. 
\end{proof}

We mention two immediate corollaries that will be useful in the sequel.
The first specializes the result to the case $m=1$:
\begin{corollary}\label{cor:classification_m1}
	\label{com:m1} 
Assume that the SHGH conjecture holds.
Let $\calL$ be an effective $n$-point linear system 
with $g(\calL) \geq 2$ and $m(\calL)=1$ 
(i.e. the biadjoint system $|C+2K_n|$ is empty,
and $\alpha=\alpha(\calL) = \dim(|C+K_n|)$).
Then $\calL$ is Cremona equivalent to 
one of the following $n$-Cremona minimal systems:
\\
\begin{inparaenum}[(i)]
\item $|3+\alpha;1+ \alpha,1^{n-1}|$;\\
\item $|5; 1^n |$ if $\alpha=5$;\\
\item $(3+\frac{\alpha} 2;
\frac {\alpha} 2,1^{n-1})$ if $\alpha\ge 2$ even;\\
\item $|3+\frac {\alpha+1} 2;\frac {\alpha+1} 2,2,1^{n-2}|$ if $\alpha\ge 3$ odd.
\end{inparaenum}\\
In the first case $|C'|=|C+K_n|$ is trivial (if $\alpha=0$) 
or composed with a pencil (if $\alpha>0$); 
in the last three, $C'^2>0$.
\end{corollary}

The second corollary is our stepping stone to prove the finiteness of Cremona minimal
linear systems of fixed genus:

\begin{corollary}\label{cor:im} 
Assume that the SHGH conjecture holds.
Fix an integer $g \geq 2$.
Consider all effective $n$-point linear systems $\mathcal L$ on $X_n$ 
with $g(\calL)=g$. 
Suppose that the set of all integers 
$\{m(\calL), \alpha(\calL): \text {$\mathcal L$ as above}\}$  
is bounded  by a constant depending on $g$.  
Then the number of points $n$ is also bounded,
and there are only finitely many 
numerical characters $|a;b_1,\ldots,b_n|$
for $n$-Cremona minimal linear systems $\mathcal L$ with $g(\calL)=g$.
\end{corollary}

\begin{proof}
Using the classification given in Theorem \ref{thm:mm},
we see that if $n$, $m$, and $\alpha$ are bounded,
then there are finitely many $n$-Cremona minimal linear systems as claimed.
We claim that if $m$ and $\alpha$ are bounded, then $n$ is also.

This is clear, after observing that the degree $d$ of the linear systems are bounded by functions of $m$ and $\alpha$.
Once $d$ is bounded, then $n$ cannot be more than $d(d+3)/2$,
since all multiplicities are at least $1$, and $n > d(d+3)/2$
would cause the linear system to be empty.
\end{proof}

\section{Finiteness for fixed genus}\label{sec:finiteness}
In this section we prove our first main result.

\begin{proof}[Proof of Theorem \ref{thm:main}] 
By the results in \cite {CDM}, if $m\leq 2$, there are only finitely many possibilities for Cremona minimal linear systems of a given genus $g(\calL)$. On the other hand if  $\calL^2>3g-3$ then $m\leq 2$ (see \cite [Lemma 3.2] {CDM}).  So we will assume from now on $\calL^2\leq 3g-3$ and $m\geq 3$.

Taking into account  Corollary \ref {cor:im}, to prove this theorem it suffices to show that the invariants $m=m(\calL)$ and $\alpha=\alpha(\calL)$ are bounded above by a function of $g$.  

First we examine the cases $n\leq 9$, in which $-K_n$ is effective. Actually, for $n\leq 8$, $-K_n$ is ample, so that for any curve $C$ on $X_n$ one has $C\cdot K_n<0$. For $n=9$ one has $(K_9)^2=0$. So, by the index theorem, for any  irreducible curve $C$  on $X_9$ such that $C^2>0$, one has $C\cdot K_9<0$. In any event, if $n\leq 9$ and $C\in \calL$, we have $C\cdot K_n<0$, i.e., $C^2>2g-2$. Now $C\cdot (C+tK_n)=C^2-t(C^2-2g+2)$ so that $|C+tK_n|$ is empty as soon as 
$$
t>\frac {C^2}{C^2-2g+2}=1+\frac {2g-2}{C^2-2g+2}.
$$
But
$$
1+\frac {2g-2}{C^2-2g+2}\leq 2g-1,
$$
and therefore $m\leq 2g-1$. Moreover $|C+mK_n|$ is a sublinear system of $|C|$, hence $\alpha\leq \dim (|C|)=C^2-g+1\leq 2g-2$ because we are assuming now $C^2\leq 3g-3$. In conclusion, if $n\leq 9$, $m$ and $\alpha$ are bounded, and we are done.

Next, we assume $n\geq 10$.
Consider an adjoint linear system $|C+tK_n|$, with $C\in \calL$ and suppose it is effective. Then
\begin{equation}\label{eq:dim}
\dim(|C+tK_n|)={\rm virtdim}(|C+tK_n|)+ h^1(X_n, \mathcal O_{X_n}(C+tK_n)),
\end{equation}
where ${\rm virtdim}(|C+tK_n|)$ is easily computed to be
\begin{equation}\label{eq:virt}
{\rm virtdim}(|C+tK_n|)=C^2-g+1+t\Big (2g-2-C^2-\frac {9-n}2\Big )+ t^2\frac {9-n}2.
\end{equation}
If there is a non--zero $h^1(X_n, \mathcal O_{X_n}(C+tK_n))$, this means that there are $h\leq n$ disjoint $(-1)$--curves $A_1,\ldots, A_h$ on $X_n$ such that  $0>(C+tK_n)\cdot A_i$ for $1\leq i\leq h$, so that $e_i:=C\cdot A_i<t$ (see Lemma \ref {lem:a}). Then, since we are assuming that SHGH Conjecture holds,  one has
$$
h^1(X_n, \mathcal O_{X_n}(C+tK_n))=\sum_{i=1}^h h^1(A_i, \mathcal O_{A_i}(C+tK_n))=\sum_{i=1}^h (t-e_i-1)\leq n(t-2)< tn.
$$
So, by \eqref {eq:dim} and \eqref {eq:virt} we have
$$
0\leq \dim(|C+tK_n|)< C^2-g+1 +t\Big (2g-2-C^2-\frac {9-3n}2\Big )+ t^2\frac {9-n}2,
$$
and note that the rightmost term is a quadratic in $t$ with negative top coefficient. So  that $t$ must be less than the largest root of that quadratic.

Taking into account that $g-1\leq C^2\leq 3g-3$, we have 
$$
0\leq \dim(|C+tK_n|)< 2g-2 +t\Big (g-1-\frac {9-3n}2\Big )+ t^2\frac {9-n}2,
$$
so
$$
2-2g -t\Big (g-1-\frac {9-3n}2\Big )+ t^2\frac {n-9}2> 0.
$$

Let us set 
$$
\Delta=\Big(g-1-\frac {9-3n}2\Big )^2+4(n-9)(g-1)
$$
so that
$$
t\leq \frac {g-1-\frac {9-3n}2+\sqrt \Delta}{n-9}
$$
and therefore also
$$
m\leq \frac {g-1-\frac {9-3n}2+\sqrt \Delta}{n-9}.
$$
Now an easy computation shows that
$$
\frac {g-1-\frac {9-3n}2+\sqrt \Delta}{n-9}\leq \varphi(g):=g+11+\sqrt{(g+11)^2+4(g-1)}
$$
so that $m\leq \varphi(g)$ is bounded by a function of $g$. 

Similarly, we have
$$
\alpha=\dim(|C+mK_n|)< 2g-2+m\Big (g-1-\frac {9-3n}2\Big )+ m^2\frac {9-n}2.
$$
Since $m\geq 3$, we have $m^3\geq 3m$ and therefore
\begin{align*}
\alpha=&\dim(|C+mK_n|)< 2g-2+m\Big (g-1-\frac {9-3n}2\Big )+ 3m\frac {9-n}2=\\
=&(m+2)(g-1)+ 9m,
\end{align*}
which proves that $\alpha$ is also bounded by a function of $g$, as desired. \end{proof}

\section{Systems of minimal selfintersection}\label{sec:bound_selfintersection}

In this section we prove Theorem \ref{thm:bound_selfintersection}, 
by explicitly determining the minimal selfintersection of a linear system 
satisfying assumptions 
\eqref{hypothesis:SHGH}, \eqref{hypothesis:g}, 
and \eqref{hypothesis:no-less-than-n-points},
given the dimension $r$ and the number of points $n$.
We then may take advantage of the consequences 
\eqref{hypothesis:negcurves}--\eqref{hypothesis:adjoint_nofix}.
We will not repeat these assumptions and consequences
in the supporting Lemma and Proposition statements,
but we will in the culminating Theorems.

It will be useful to collect below the useful classifications
of systems of genus at most two;
these are classical results, 
for which one may consult Section 10 of \cite{CCi} for a modern treatment:

\begin{proposition}\label{prop:classg=012}
Assume that the SHGH Conjecture holds, 
and $|C|$ is a nef effective Cremona minimal system on $X_n$ 
with irreducible general member of genus $g$. Then:
\begin{itemize}
\item[(0)] If $g=0$, then $|C|$ is one of the systems 
$|2;0^n|$,
$|t;t-1,0^{n-1}|$ for $t \geq 1$, or
$|t;t-1,1,0^{n-2}|$ for $t \geq 1$.
\item[(1)] If $g=1$, then $|C|$ is one of the systems
$|3;1^k,0^{n-k}|$ for $0\leq k \leq 9$ ($r=9-k$); or
$|4;2,2,0^{n-2}|$ ($r=8$).
\item[(2)] If $g=2$, and $r \geq 1$, then $|C|$ is one of the systems
$|4;2,1^k,0^{n-k-1}|$, $0 \leq k \leq 10$, $n \geq k+1$ ($r=11-k$);
$|6;2^8,1^k, 0^{n-k-8}|$, $0 \leq k \leq 2$, $n \geq k+8$ ($r=3-k$);
\end{itemize}
\end{proposition}

\subsection{The hyperelliptic case}\label{sec:hyper}
	
In this subsection we will consider the case in which 
either $|C'|=|C+K_n|$ is a pencil (hence $g=2$), 
or the general curve in $|C'|$ is reducible,  
which by consequence \eqref{hypothesis:adjoint_nofix}
implies that $|C'|$ is composed with a pencil $|F|$. 
In this case the curves in $|C|$ are hyperelliptic and $C\cdot F=2$.
Since $\dim(|C'|) = g-1 \geq 1$, in this case $C'$ is nef,
so that by Proposition \ref{prop:fin} we see that
if it is the case that $C'|$ is composed with a pencil, 
equal to $|\beta F|$ with $\beta = g-1 \geq 2$, then
$F^2=0$ (and $|F|$ is Cremona equivalent to $|1;1,0^{n-1}|$, with genus zero). 
	
\begin{lemma}\label{lem:one} 
	$H^1(F,\mathcal O_F(C))=0$.
\end{lemma}
	
\begin{proof} 
Consider the exact sequence
$$
		0 \longrightarrow \mathcal O_{X_n}(C-F)
		\longrightarrow \mathcal O_{X_n}(C)
		\longrightarrow \mathcal O_F(C)
		\longrightarrow 0.
$$
By consequence (G), we have $H^1(X_n, O_{X_n}(C))=0$,
so that using the long exact sequence
it suffices to prove that $H^2(X_n,\mathcal O_{X_n}(C-F))=0$.
By Serre duality, this is equivalent to having 
$H^0(X_n, \mathcal O_{X_n}(K_n+F-C)) = 0$,
and so we must show that $K_n+F-C$ is not effective.
		
Since $F$ is a component of $C'\equiv K_n+C$,
we have that $K_n+F-C \leq K_n+C'-C \equiv 2K_n$,
which is never effective, and so neither is $K_n+F-C$.
\end{proof}

Let $p$ be the arithmetic genus of $F$.
By the Lemma \ref{lem:one}, since $C\cdot F = 2$,
we have $h^0(\calO_F(C)) = 2+1-p = 3-p$,
and is an effective divisor of degree $2$, so $p\leq 2$.

\begin{lemma}\label{lem:two} 
If $p=2$, then $r=0$.
If $p \geq 1$, then $g=2$; hence $\beta=1$ and $|C'|=|F|$.
\end{lemma}
	
\begin{proof}
Suppose that $p=2$.
Then $h^0(F,\mathcal O_F(C))=1$, and hence if $r \geq 1$ (so that $C$ moves),
$|C|$ must have two base points on the general curve $F$.
Since $F$ moves, these points must be base points for $|F|$ also.
However then every element of $|F|$ meets the general $C$
in these two points, contradicting that $|F|$ cuts out the
hyperelliptic pencil on the general $C$.
We conclude that if $p=2$ then $r=0$.

The second statement follows from the observation that if $\beta \geq 2$ then $p=0$;
hence if $p \geq 1$ then $\beta = 1$, forcing $g=2$.
\end{proof}
	
\begin{theorem}\label{thm:hyperelliptic}
Assume that the SHGH conjecture holds.
Let $\mathcal L$ be a linear system on $X_n$ 
satisfying \eqref{hypothesis:g} and \eqref{hypothesis:no-less-than-n-points}.
Let $g=g(\calL)$ and $r=\dim(\calL)$.
Suppose that the adjoint linear system $|C+K_n|$
is a pencil or is composed with a pencil,
(which is forced if $g=2$). 
Then $\mathcal L$ is Cremona equivalent to 
one of the following Cremona minimal linear systems:
\begin{enumerate}[(i)]
\item $|6;2^8, 1^{n-8}|$ with $8\leq n \le 11$; in this case 
	\begin{itemize}
	\item $n+r=11$;
	\item $g=2$;
	\item $C^2 = 12-n$.
	\end{itemize}
\item $|g+2;g,1^{n-1}|$ with $1\le n \le 3g+6$; in this case 
	\begin{itemize}
	\item $n+r=3g+6$ is divisible by $3$;
	\item $g = (n+r)/3 - 2$;
	\item $C^2 = 4g-n+5 = (n+4r)/3 - 3$.
	\end{itemize}
\item $|9;3^8,2^2|$, for which $n=10$, $r=0$, $g=2$, and $C^2 = 1$.
\end{enumerate}	
\end{theorem}
	
\begin{proof}
We consider the three cases, 
depending on the genus $p=g(F)$.
We let $m=m(\calL)$ and $\alpha=\alpha(\calL)$ as defined above. 

If $p=2$, then $r=0$, $g=2$, and $|C'|=|F|$ by Lemma \ref{lem:two}.
We also have that $|F|$ has dimension one (it is a pencil)
so that Proposition \ref{prop:classg=012}(2) applies to $|F|$, since $F$ is nef;
we conclude that $|F|$ must be one of the dimension one systems 
$|4;2,1^{11},0^{n-12}|$ or $|6;2^8,1^2,0^{n-10}|$.
Hence $|C=F-K|$ must be either
$|7;3, 2^{11},1^{n-12}|$ or 
$|9;3^8,2^2,1^{n-10}|$.
The degree $7$ system is not effective, and we conclude that $|C|$ is this degree $9$ system, with only $n=10$ giving $r=0$.
		
If $p=1$, we again have $g=2$ and $|C'|=|F|$.
Then, up to a Cremona transformation, 
$|F|$ is the pencil $|3;1^8, 0^{n-8}|$ 
(see \cite[Sect. 5.4] {Ca})
and therefore $|C|=|F-K_n|$ 
is Cremona equivalent to 
$|6;2^8, 1^{n-8}|$, with $n\geq 8$.  
This system has dimension $r=27-24-n+8=11-n$, so $n\leq 11$.
		
If $p=0$, that is, $|C'|$ is composed with a pencil of rational curves, 
then $|C'+K|=|C+2K|$ is empty and therefore $m(\calL)=1$.
Therefore Corollary \ref{cor:classification_m1} applies,
and we must be in the first case (i) there, since $|C'|$ is composed with a pencil. 
Since $\alpha = g-1$,
this gives that $|C|$ is Cremona equivalent to $|g+2;g,1^{n-1}|$.
One computes then that $C^2 = 4g-n+5$ and $r = 3g-n+6$, 
completing the proof.
\end{proof}

Note that, in the hyperelliptic systems classified in this subsection, 
the self-intersection is determined
by $n$ and $r$.
We will see that this is not true in the non-hyperelliptic case.

\subsection{The self-intersection of non-hyperelliptic systems}\label{sec:si-nonh}
The question of the minimum self-intersection possible,
for a given $n$ and $r$, is more complex for the non-hyperelliptic case.
In this case, as we will see below, the self-intersection
depends on $n$, $r$, and the genus of the adjoint system.

Although this was not needed in the hyperelliptic case analysis of the previous subsection, for the non-hyperelliptic cases we will assume that we are not in the weak del Pezzo situation: we will assume that $n\geq 10$.

Consider the case in which 
$|C|$ is an effective complete linear system on $X_n$ 
of dimension $r\geq 0$ with $g\geq 3$,
and its adjoint system $|C'|$ is not composed with a pencil, 
so that by \eqref{hypothesis:adjoint_nofix}, 
the general curve of $|C'|$ is irreducible.
We may expect that in this case the general curve in $|C|$ is not hyperelliptic,
though this is not really guaranteed. 
However, as we will a posteriori see, if $C^2$ is minimal (or close to minimal), 
in this case the general curve in $|C|$ is indeed not hyperelliptic. 
Therefore we will call this the \emph{non--hyperelliptic} case.

We denote by $g'$ the arithmetic genus of $C'$. 

\begin{proposition}\label{prop:max} In the above set--up, we have
\[
	C^2=\frac {n+3r+g'}{2} - 5.
\]
\end{proposition}

\begin{proof} Consider the exact sequence
	\[
	0 \longrightarrow \mathcal O_{X_n} (-K_n)
	\longrightarrow  \mathcal O_{X_n} (C)
	\longrightarrow  \mathcal O_{C'} (C)
	\longrightarrow 0.
	\]
Since $h^1(X_n, O_{X_n} (C))=0$ 
(consequence \eqref{hypothesis:cohomology})
and 
$$
h^2(X_n, \mathcal O_{X_n} (-K_n))=h^0(X_n,\mathcal O_{X_n} (2K_n))=0,
$$ 
the long exact sequence implies that $h^1(C', \mathcal O_{C'} (C))=0$. 
Since  $h^0(X_n, \mathcal O_{X_n} (-K_n))=0$
(here we use $n \geq 10$),
we have
\[
	h^0(C', \mathcal O_{C'} (C))=r+1+h^1(X_n, \mathcal O_{X_n} (-K_n)).
\]
	As $\chi(-K_n)=K_n^2+1=10-n$, we find
\[
	h^0(C', \mathcal O_{C'} (C))=r+1-(10-n)=r+n-9.
\]
Since, as we saw, $h^1(C, \mathcal O_{C} (C))=0$ 
and $h^0(C, \mathcal O_{C} (C))=r$, 
we deduce that $r=C^2-g+1$, i.e., $g=C^2-r+1$. 
One has
	\[
	C\cdot C'=2g-2=2C^2-2r.
	\]
	Hence
	\[
	r+n-9=h^0(C', \mathcal O_{C'} (C))=2C^2-2r-g'+1
	\]
	and therefore
	\[
	2C^2=3r+n-10+g'
	\]
	as claimed. 
\end{proof}

Because of this Proposition,
the analysis splits in two, according to the parity of $n+r$,
which is that of $g'$.
We take up these two cases in the next two subsections,
where we may assume that $g\geq 3$.

In seeking the minimum self-intersection for such linear systems,
we clearly seek such systems with minimum $g'$.
If we write $n+r = 2h+\epsilon$, with $\epsilon \in \{0,1\}$,
by Proposition \ref{prop:max} we have
\begin{equation}\label{eq:self}
	C^2=h+r-5+\frac {g'+\epsilon}2.
\end{equation}
Then $g'=\epsilon$ is the minimum possible
(that gives integer values for $C^2$).

We note that, using the sequence
\begin{equation}\label{eq:ses1}
0\longrightarrow \mathcal O_{X_n}\longrightarrow \mathcal O_{X_n}(C)\longrightarrow \mathcal O_C(C)\longrightarrow 0,
\end{equation}
we see that $h^0(\calO_C(C)) = r$, and has degree $C^2$;
it also has no $H^1$, using \eqref{hypothesis:cohomology}.
hence by Riemann-Roch, we have $r = C^2+1-g$,
or $C^2-r+1 = g \geq 2$.
Then $C^2 \geq r+1$,
which implies,  by Proposition \ref{prop:max}, 
that $(n+r+g')/2 \geq 6$, 
or $h + (g'+\epsilon)/2 \geq 6$.

In any case we see the following immediate conclusion, from \eqref{eq:self}:
\[
C^2 \geq \begin{cases}
h+r-5 & \text{ if } n+r = 2h \text{ is even} \\
h+r-4 & \text{ if } n+r = 2h+1 \text{ is odd}
\end{cases}
\] 

\subsection{Non-hyperelliptic systems with $n+r$ even}
\label{sec:nonh-n+r=even}

In this subsection we will prove that in the case that $n+r$ is even,
the bound given above is almost always achieved,
and we will classify the systems that achieve the bound
(up to Cremona equivalence).
Any system that achieves the bound must have $g'=0$,
and we start by analyzing such systems.


\begin{lemma}\label{g'=0}
In the non-hyperelliptic case with $g \geq 3$, 
if $g'=0$, $n \geq 10$, 
and $|C|$ is $n$-Cremona minimal,
then $\epsilon=0$, $n+r=2h$ is even,
and $|C|$ is equal to one of the following:
\begin{itemize}
	\item[(a)] $|5;1^n|$;
	here $r=20-n$, $h=10$, $g=6$, and $C^2 = 25-n$, 
	for $10 \leq n \leq 20$.
	\item[(b)] $|t+3;t, 1^{n-1}|$ 
	for $t \geq 1$ and $10 \leq n \leq 4t+10$;
	here $h$ is odd, equal to $2t+5$, with $r=4t-n+10$, $g=2t+1$, 
	and $C^2=6t-n+10$.
	\item[(c)] $|t+3;t,2,1^{n-2}|$
	for $t \geq 2$ and $10 \leq n \leq 4t+8$; 
	here $h$ is even, equal to $2t+4$, 
	with $r=4t-n+8$, $g = 2t$, and $C^2 = 6t-n+7$.
\end{itemize}
\end{lemma}

\begin{proof}
The non-hyperelliptic systems with rational adjoint system are classified,
see \cite[Sect. 5.3]{Ca}.  The result is that
$|C'|$ must be (up to a Cremona transformation),
one of the systems of Proposition \ref{prop:classg=012}(0).
Therefore $|C| = |C'-K|$ is Cremona equivalent to either
$|5;1^n|$;
$|t+3;t, 1^{n-1}|$, or
$|t+3;t,2,1^{n-2}|$ as claimed,
with the stated values therefore for $r$, $g$, $h$, and $C^2$ in each case.
(The restriction in (c) that $t \geq 2$ is because 
if $t=1$ then the system is the hyperelliptic system $|4;2,1^{9}|$, 
which has $g=2$, and we are assuming $g \geq 3$.) 
\end{proof}

Lemma \ref{g'=0} now allows us to state the main result here.

\begin{theorem}\label{thm:irreduciblen+reven}
Suppose that the SHGH Conjecture holds,
as well as assumptions \eqref{hypothesis:g}, 
and \eqref{hypothesis:no-less-than-n-points}
for an effective linear system $|C|$.
Assume further that $n \geq 10$, $g \geq 2$, and $n+r = 2h$ is even.
Then:
\begin{itemize}

\item[(A)] If we are in the non-hyperelliptic case,
so that the adjoint linear system $|C'|$ has irreducible general member, 
then $n+r \geq 14$, 
and the minimum self-intersection $C^2$ 
with these assumptions is $h+r-5$.
An $n$-Cremona-minimal linear system achieves this minimum 
if and only if 
it is one of the systems listed in Lemma \ref{g'=0}.
All values of $(n,r)$ with $n+r \geq 14$ occur on this list,
i.e., for any $n \geq 10$ and $r \geq 0$ with $n+r=2h\geq 14$ even,
$n$-Cremona minimal systems $|C|$ with $C^2 = h+r-5$ exist 
and all are listed above.
\item[(B)] 
If $n+r=10$ or $12$, then $g=2$, with the following cases: 
\begin{itemize}
\item[(10,0)] If $n=10$ and $r=0$, the minimum $C^2$ for any system (hyperelliptic or not) is $1$; 
the only Cremona-minimal system with $C^2=1$ is the hyperelliptic system $|9;3^8,2^2|$. 
\item[(12,0)] If $n=12$ and $r=0$, the minimum $C^2$ for any system (hyperelliptic or not) is $1$;
the only Cremona-minimal system with $C^2=1$ 
is the hyperelliptic system $|4;2,1^{11}|$.
\item[(11,1)] If $n=11$ and $r=1$, the minimum $C^2$ for any system
(hyperelliptic or not) is $2$;
the only Cremona-minimal system with $C^2=2$ 
is the hyperelliptic system $|4;2,1^{10}|$.
\item[(10,2)] If $n=10$ and $r=2$, the minimum $C^2$ for any system
(hyperelliptic or not) is $3$;
the only Cremona-minimal system with $C^2=3$ 
is the hyperelliptic system $|4;2,1^9|$.
\end{itemize}
\end{itemize}
\end{theorem}

\begin{proof}
To prove (A), we need only address the final statement,
that all such $(n,r)$ with $n+r \geq 14$ are possible.
If we write $n+r=2h$ and $h$ is odd, then since $n+r \geq 14$,
we must have $h \geq 7$.
Hence we may set $t = (h-5)/2$ for some $t \geq 1$.
In that case the system $|t+3;t,1^{n-1}|$ of (b) exists with the given $n$ and $r$:
we have $r = 4t-n+10 = 2(h-5)-n+10 = 2h-n$ as required.
	
If $h$ is even, we must have $h \geq 8$;
in this case set $t = (h-4)/2$.
Then the system $|t+3;t,2,1^{n-2}|$ of (c) exists with the given $n$ and $r$:
$r=4t-n+8 = 2(h-4)-n+8 = 2h-n$ as required. 

To prove (B)(10,0), we then have $h=5$,
and the bound above is that $C^2 \geq 0$.
If $C^2=0$, then $g'=0$, so Lemma \ref{g'=0} applies;
none of those systems have $r=0$ though.
If $C^2=1$, then $g'=2$, and $r'=\dim(|C'|) \geq g-1 \geq 2$,
so the classification of genus two systems 
from Proposition \ref{prop:classg=012}(2)
gives that $|C'|$ must be one of the systems listed there,
and we conclude that $|C|$ is one of the systems
$|7;3,2^k,1^{9-k}|$ ($r=35-6-3k-(9-k)=20-2k$, 
$|8;4,3,1^{8}|$ ($r=44-10-6-8=20$), 
$|9;3^8,2^k,1^{2-k}|$ ($r=54-48-3k-(2-k)=4-2k$,
$|12;4^8,3^2|$ ($r=90-80-12<0$).
The system of degree $7$ does not have $r=0$ for any $k\leq 9$;
the system of degree $8$ has $r > 0$,
and the system of degree $12$ does not exist
The system of degree $9$ with $k=2$ is the only possibility:
$|9;3^8,2^2|$ and this does have $C^2 = 81-72-8 = 1$
(but is hyperelliptic as noted in Theorem \ref{thm:hyperelliptic}).

%

For the case $(n,r)=(12,0)$, we have $h=6$,
so the lower bound for $C^2$ is $1$.
If the system is not hyperelliptic, with $C^2 = 1$,
then we must have $g'=0$,
and again Lemma \ref{g'=0} applies, with none of the systems having $r=0$.
Hence $C^2 \geq 2$ for any non-hyperelliptic system.

The hyperelliptic system $|4;2,1^{11}|$ has $C^2=1$ and is therefore the unique system with this minimal self-intersection, for $(n,r)=(12,0)$.

For the case (B)(11,1), we have again $h=6$, so that the bound
for a non-hyperelliptic system is that $C^2 \geq 2$.
If $C^2 = 2$, then $1=r=C^2+1-g = 3-g$ so $g= 2$
and the system cannot be non-hyperelliptic.
Therefore we conclude that the minimum self-intersection
for any such system will be $2$,
and occurs only with the hyperelliptic system $|4;2,1^{10}|$.

Finally for the case (B)(10,2), again we have $h=6$,
so the bound for a non-hyperelliptic system gives $C^2 \geq 3$.
However if $C^2 = 3$,
then $2=r=C^2+1-g = 4-g$ and again $g=2$, forcing the hyperelliptic system
$|4;2,1^9|$.

%
\end{proof}

\begin{remark}\label{rem:notunique} In the case $n+r$ even, examined above, for $n\geq 19$ the linear system $|C|$ of dimension $r$ with $|C'|$ not composed with a pencil and with minimal self intersection is unique up to Cremona transformation. For $10 \leq n\leq 18$ the system is not unique. Indeed we have the system $|5;1^n|$ of case (a) of Lemma \ref{g'=0}, 
	with $r=20-n$ and self intersection $25-n$. 
	But we have also the system $|6;3,2, 1^{n-2}|$ of case (c) of that Lemma, 
	($t=3$ with this system)
	with the same dimension and self intersection
	and the two are Cremona minimal and not Cremona equivalent.  
\end{remark}

\subsection{Non-hyperelliptic systems with $n+r$ odd}\label{sec:nonhn+rodd}
If $n+r$ is odd, by Proposition \ref{prop:max},
we will achieve the minimum self-intersection for the
non-hyperelliptic case
if we can arrange that $g'=1$;
using the notation of the prior section, we will then have $\epsilon=1$,
so that we may write $n+r = 2h+1$.
We would then be seeking systems $|C|$ 
with $C^2 = h+r-4$.

We can try to analyze this case as in the $n+r$ even case, 
where we looked for $g'=0$ systems.  
In that case we were able to find
systems for every relevant pair of integers $(n,r)$ 
with $n+r \geq 14$.
However we will see that the analysis ends up being quite different: 
with $n+r$ odd there are only finitely many systems with $g'=1$.

\begin{lemma}\label{g'=1}
We assume $n\geq 10$ and $g \geq 3$ as above.
In the non-hyperelliptic case, if $g'=1$, 
and $|C|$ is $n$-Cremona minimal,
then $\epsilon=1$, $n+r=2h+1$ is odd,
and $|C|$ is equal to one of the following: 
\begin{itemize}
	\item[(a)] $|6;2^m,1^{n-m}|$, for some $m \in \{0,\ldots,7\}$;
		here $r=27-2m-n$, $h=13-m$, $g=10-m$, 
		and $C^2 = 36-3m-n$,
		with $n \leq 27-2m$.
	\item[(b)] $|7;3^2, 1^{n-2}|$;
		here $r=25-n$, $h=12$, $g=9$, 
		and $C^2=33-n$,
		with $n \leq 25$.
\end{itemize}
\end{lemma}

\begin{proof} If $g'=1$ one has $\epsilon=1$ by \eqref {eq:self}. 
In this case of $g'=1$, as we noted in the Introduction, 
$|C'|$ must be (up to a Cremona transformation),
a linear system of the form
$|3;1^m|$, with $0\leq m\leq 9$, or $|4;2^2|$.
Therefore $|C| = |C'-K|$ is Cremona equivalent to either
$|6;2^m,1^{n-m}|$, or
$|7;3^2, 1^{n-2}|$, as claimed,
with the stated values therefore for $r$, $g$, $h$, and $C^2$ in each case.
In case (a) we must have $m \leq 7$ since $g \geq 3$.
The bounds on $n$ are a result of $r \geq 0$ in both cases.
\end{proof}




\begin{remark}\label{rem:notu} Also in this case with $n+r$ odd, 
	the linear system $|C|$ of dimension $r$ 
	with $|C'|$ not composed with a pencil and with minimal self intersection 
	is not unique up to Cremona transformation for $g=9$. 
	In  fact we have the two systems 
	$|6;2, 1^{n-1}|$ and $(7, 3^2, 1^{n-2})$ for 
$10\leq n\leq 25$,
	both having $r=25-n$ and minimal self intersection $33-n$, 
	both Cremona minimal and not Cremona equivalent. 
\end{remark}

We see then that for large $n$ and $r$ with $n+r$ odd,
we cannot hope to achieve the value of $h+r-4$
for $C^2$; this would require $g'=1$ and as seen above that leads to
systems with bounded $n$ and $r$.

We could try to analyze the $g'=3$ cases
but it turns out that this will also not produce all examples;
indeed, in this case of $n+r$ odd,
we need to allow $g'$ to increase without bound as $n$ and $r$ increase.

To see this we require a more refined lower bound for $C^2$
than Proposition \ref{prop:max} provides.
For this purpose we recall the following result
originally due to Castelnuovo and Enriques
(see \cite {CR} for a modern reference):

\begin{theorem}\label{thm:CE} 
	Let $X$ be a smooth, irreducible, projective surface. 
	Let $D$ be an irreducible curve on $X$. 
	Set $\rho := \dim(|D|)$ and let $\gamma$ be the arithmetic genus of $D$.
	Assume that there is no smooth rational curve $F$ on $X$
	such that $F^ 2 = 0$ and $D \cdot F = 1$.  
	Then $\rho \leq 3\gamma+5+\eta$, 
	where $\eta=1$
	if $\gamma=1$ otherwise $\eta=0$.
\end{theorem}

We apply this theorem to the adjoint linear system $|C'|$,
in order to prove the following.

\begin{proposition}\label{prop:CEboundC'}
	Fix $n \geq 10$, $r \geq 0$, with $n+r=2h+1$ being odd.
	Suppose that the general member of the adjoint linear system is irreducible
	of genus $g' \geq 3$.
	Then
	\begin{equation}\label{eq:erot}
		C^2 \geq r+ \left\lceil{\frac{6h-7}{5}}\right\rceil - 5.
	\end{equation}
\end{proposition}

\begin{proof}
	Since $n+r$ is odd we must have $g'$ odd by Proposition \ref{prop:max};
	let us write $g' = 2a-1$ with $a \geq 2$.
	Then $C^2 = (n+r+g')/2 + r - 5 = h+a+r-5$.
	Riemann-Roch for $\calO_C(C)$ implies that $r-1 = C^2 - g$,
	which  implies  $g = h+a-4$.
	However we also have $g = \dim(|C'|)+1$.
	We can apply Theorem \ref{thm:CE} with $D=C'$,
	since the surface $X_n$ is rational and $C'$ is not,
	so that there is certainly no smooth rational curve $F$ with $C'\cdot F = 1$.
	In that case we have (in the notation of the Theorem)
	$\rho = g-1$, $\gamma = g'$, and $\eta = 0$,
	and we conclude that
	$g-1 \leq 3g' + 5$, or $h+a-5 \leq 3(2a-1) + 5$,
	i.e., $h \leq 5a+7$, or
	\[
	a \geq \left\lceil{\frac{h-7}{5}}\right\rceil.
	\]
	Hence
	\[
	C^2 = h+a+r-5 \geq h+r-5 + \left\lceil{\frac{h-7}{5}}\right\rceil
	= r+ \left\lceil{\frac{6h-7}{5}}\right\rceil - 5
	\]
	as claimed.
\end{proof}

\begin{remark}\label{rem:for} 
	Keep the notation of Proposition \ref {prop:CEboundC'}.
	If we set 
	$h-7=5b-m$, with $0\leq m\leq 4$, then we have
	$$
	\left \lceil {\frac {h-7}5} \right \rceil=b
	$$
	and we have equality in \eqref {eq:erot} if and only if $a (=(g'+1)/2)=b$. 
\end{remark}

We see from the above argument that we cannot hope to bound $g'$
as we could for the $n+r$ even case, since the quantity $a$ grows with $h$.

We now show that the bound for $C^2$ given in the above Proposition is sharp.

\begin{lemma}\label{lem:2b+4system}
Consider the linear system $|C|=|2b+4;2b,2^m,1^k|$
with $b \geq 1$, $0\leq m \leq 4$, and $9-m \leq k \leq 10b-3m+14$.
Then:
\begin{align*}
	n &= m+k+1 \\
	r &= 10b-3m-k+14 \\
	n+r &= 10b-2m+15, \text{ so that } h = 5b-m+7 \\
	g &=6b-m+3  \\
	g' &=2b-1,  \text{ so that } a=b\\
	C^2 &= 16b-4m-k+16
\end{align*}
and this system achieves the bound of Proposition \ref{prop:CEboundC'}.
\end{lemma}

\begin{proof}
The inequalities on $k$ ensure that $n \geq 10$ and $r \geq 0$.
We leave to the reader the calculation of these quantities for the system;
we simply note that, with these parameters, we have
$6h-7 = 30b-6m+35 = 5(6b-m+7) - m$ so that
\[
	\left\lceil{\frac{6h-7}{5}}\right\rceil = 6b-m+7
\]
for all values of $m \in \{0,\ldots,4\}$.
Hence the lower bound is $r+(6b-m+7) - 5 = 16b-4m-k + 16$,
which is the value for $C^2$. 
\end{proof}

\begin{theorem}\label{thm:irreduciblen+rodd}
Suppose that the SHGH Conjecture holds,
as well as assumptions \eqref{hypothesis:g}, 
and \eqref{hypothesis:no-less-than-n-points}
for an effective linear system $|C|$.
Assume further that $n \geq 10$, $g \geq 3$, $n+r = 2h+1$ is odd,
and the adjoint linear system $|C'|$ has irreducible general member of genus $g'$.
Then:
\begin{itemize}
\item[(A)] If $n+r \geq 13$, the minimum self-intersection $C^2$ 
with these characters is the bound of \eqref{eq:erot}.
An $n$-Cremona-minimal linear system with these characters 
has minimal self-intersection if and only if it is one of the following:
\begin{itemize}
	\item[(a)] $|6;2^\ell,1^{n-\ell}|$, for some $\ell \in \{0,\ldots,7\}$;
		here $r=27-2\ell-n$, $h=13-\ell$, $g=10-\ell$, 
		and $C^2 = 36-3\ell-n$,
		with $n \leq 27-2\ell$.
	\item[(b)] $|7;3^2, 1^{n-2}|$;
		here $r=25-n$, $h=12$, $g=9$, 
		and $C^2=33-n$,
		with $n \leq 25$.
	\item[(c)] $|7;2^\ell,1^k|$ with $\ell+k=n$ and $0\leq \ell \leq 4$;
		here $r =35-3\ell-k$ and $C^2 = 49-4\ell-k$.
	\item[(d)] $|8;3,2^\ell,1^k|$ with $n=\ell+k+1$ and $0\leq \ell\leq 1$;
		here $r=38-3\ell-k$ and $C^2 = 55-4\ell-k$.
	\item[(e)] $|2b+4;2b,2^m,1^k|$
		with $0\leq m \leq 4$, $b \geq 1$, 
		and $9-m \leq k \leq 10b-3m+14$;
		here $n = m+k+1$, $r = 10b-3m-k+14$, 
		$n+r = 10b-2m+15$, and $C^2 = 16b-4m-k+16$.
		(These are the systems of Lemma \ref{lem:2b+4system}.)
\end{itemize}
All values of $(n,r)$ with $n\geq 10$, $r \geq 0$, $n+r \geq 13$ and odd 
occur on this list; the only missing systems with $n+r$ odd 
are the cases with $n+r=11$,  
namely the systems with $(n,r)=(10,1)$ or $(11,0)$.
\item[(B)]
\begin{itemize}
\item[(10,1)] If $(n,r)=(10,1)$, 
the non-hyperelliptic Cremona minimal systems with minimal self-intersection 
have $C^2=3$, and are either $|10;4,3^9|$ or $|12;4^8,3,2|$. 
(We note that there is the hyperelliptic system 
$|6;2^8,1^2|$ with $C^2=2$, though.)
\item[(11,0)] If $(n,r)=(11,0)$, 
the Cremona minimal system with minimal self-intersection 
have $C^2=2$, and is $|9;3^7,2^4|$.
(We note that there is the hyperelliptic system 
$|6;2^8,1^3|$ with $C^2=1$, though.)
\end{itemize}
\end{itemize}
\end{theorem}

\begin{proof}
We begin with the case (A), where $n+r \geq 13$.
Since $n+r$ is odd, we must have $g'$ odd as well,
and if $g'=1$ we have the classification of Lemma \ref{g'=1},
which gives us cases (a) and (b).
Hence we can assume that $g' \geq 3$,
and we only need to show that there are no other systems
except those of cases (c)-(e)
with these parameters, having the minimal self-intersection.
	
Consider now a linear systems $|C|$ for which the equality holds in \eqref {eq:erot}. 
By taking into account the proof of Proposition \ref {prop:CEboundC'} 
and Remark \ref {rem:for}, 
setting $b = \lceil (h-7)/5 \rceil$ and $m=5b-h+7$ as there,
we have $g'=2b-1$ and $g = h+b-4$, (hence $\dim(|C'|)=g-1=h+b-5$), so
$$
	3g'+5-\dim (|C'|)=3(2b-1)+5-(h+b-5) = 5b-h+7 = m\in \{0,\ldots, 4\}.
$$
By the Riemann--Roch theorem and since consequence \eqref{hypothesis:nobase} 
implies that  $h^1(C', \mathcal O_{C'}(C'))=0$,
one has $\dim(|C'|)=(C')^2-g'+1$, so we have
$$
	m=3g'+5-\dim (|C'|)=3g'+5-((C')^2-g'+1)=4g'+4-(C')^2,
$$
thus
$$
   4g'\leq (C')^2\leq 4g'+4. 
$$
Hence, by the Riemann--Roch theorem, if $C'\in |C+K_n|$ is a general curve, 
the linear series cut out by $|C+K_n|$ on $C'$ has no base points 
and determines a birational map to its image. 
So the map determined by the adjoint system $|C+K_n|$ 
is a birational morphism of $X_n$ to its image.  
Then if $g' \geq 7$, \cite[Thm. 1.2 and \S 3.1] {CD} implies that 
$|C'|$ is a linear system on $X_n$ of the form 
$|g'+2;g', 1^m, 0^k|=|2b+1;2b-1, 1^m, 0^k|$, 
with $n=m+k+1$ and $0\leq m\leq 4$. 
This implies that $|C|$ is one of the linear systems $|2b+4;2b,2^m,1^k|$ 
as in Lemma \ref{lem:2b+4system} and part (e) of the theorem,
for which the equality holds in \eqref {eq:erot} 
and therefore have minimal self--intersection. 
	
If $g'=3$, \cite[Thm. 1.2 and \S 3.1] {CD}  implies that for $|C'|$, 
besides the same formula as above with $b=2$
(which leads again to case (e)), 
we have another possibility. 
Namely $|C'|$ can be also of the form 
$|4;1^m,0^k|$ with $m+k=n$ and $0\leq m\leq 4$ 
(see \cite [Example 3.2]{CD}), 
so that for $|C|$ we have the  further possibility $|7; 2^m, 1^k|$ 
which leads to case (c). 
	
If $g'=5$, again \cite[Thm. 1.2 and \S 3.1] {CD}  implies that for $|C'|$, 
besides the same possibility as above with $b=3$
(which leads  to case (e)), 
we also have another possibility. 
Namely $|C'|$ can be also of the form 
$|5;2,1^m, 0^k|$ with $m+k+1=n$ and $0\leq m\leq 1$ 
(see \cite [Example 3.5]{CD}), 
so that for $|C|$ we have the  further possibility $|8; 3, 2^m, 1^k|$ 
which leads to case (d). 
	
Finally we have to prove the final statement that all $(n,r)$ are possible
as long as $n+r \geq 13$. 
	Given $n$ and $r$, we have $n+r=2h+1$. 
	Since $n+r\geq 13$, we have $h\geq 6$ 
	and if $h=6$, then either $(n,r)=(10,3)$ or  $(n,r)=(11,2)$. 
	The case $(n,r)=(10,3)$ is covered for instance, in (a) for $m=7$, $n=10$.  
	The case $(n,r)=(11,2)$ is covered in (a) for $m=7$, $n=11$.  
	If $h=7$, the possible cases are $(10,5), (11, 4), (12, 3), (13,2)$ 
that are respectively covered by the systems $|6;2^6, 1^s|$, 
with $4\leq s\leq 7$ (all in (a)). 
If $h\geq 8$, we set $h-7=5b-m$, with $b\geq 1$, $0\leq m\leq 4$; 
then the system $|2b+4;2b,2^m,1^{n-m-1}|$ has the given $n$ and $r$. 

This finishes the proof of (A); we now address (B).

We note that in the two missing cases with $n+r=11$, we have $h=5$.
If there is a system with $g\geq 3$ with the minimum self-intersection
then (in the notation of Proposition \ref{prop:CEboundC'} and Remark \ref{rem:for})
we would have $b = \lceil (h-7)/5 \rceil = 0$,
and therefore $a=0$ as well.
This would imply that $g=h+a-4 = 1$, a contradiction.
Hence in these cases a non-hyperelliptic system cannot achieve the theoretical bound.

Suppose $(n,r)=(10,1)$, so that $h=5$, with $g \geq 3$.
The bound given above is $C^2 \geq 2$, 
and we have noted above that $C^2=2$ is not possible
(since $r-1=C^2-g$, forcing $g=2$
and the system would be hyperelliptic).
Hence $C^2 \geq 3$.

If $C^2 = 3$, then $g=g' = 3$, and we compute that $C\cdot K_n=1$, $r'=2$.
The second adjoint $C''=C+2K_n$ has genus $g''=2$ and dimension $r''=2$.
In that case we have the possibilities $|C''|=|4;2,1^9|$ or $|6;2^8,1,0|$,
leading to $|C'|=|7;3,2^9|$ or $|9;3^8,2,1|$,
hence $|C|=|10;4,3^9|$ or $|12;4^8,3,2|$.

Suppose that $(n,r)=(11,0)$, so that $h=5$ with $g \geq 3$.
Again the minimal $C^2$ is larger than the bound given above (which is zero).
Since $g = C^2-r+1$, we must have $C^2 \geq 2$.

If $C^2=2$, again a computation with the second adjoint
(which will have genus one and dimension two)
leads to the system $|C|=|9;3^7,2^4|$ as claimed.
\end{proof}

\begin{remark} 
We note that the minimum $C^2$ grows like:
\begin{itemize}
	\item[(a)] $(1/3)n + (4/3)r$ in the hyperelliptic case;
	\item[(b)]  $(1/2)n+(3/2)r$ in the non-hyperelliptic case, $n+r$ even; 
	\item[(c)] $(3/5)n+(8/5)r$ in the non-hyperelliptic case, $n+r$ odd.  
\end{itemize}
\end{remark}

\section{Remarks, Applications, Observations}\label{sec:final}

In this final section we assume the hypotheses 
\eqref{hypothesis:SHGH}--\eqref{hypothesis:no-less-than-n-points}, 
which then imply that we have the consequences
\eqref{hypothesis:negcurves}--\eqref{hypothesis:adjoint_nofix}.
We will also uniformly assume that $n\geq 10$ and $r \geq 0$.
We will refer to these as the \emph{full set of hypotheses}.

\subsection{Systems with low self-intersection}

We now want to classify those systems with $C^2 \leq 5$,
whether or not they have minimal $C^2$ for the given $n$ and $r$,
which we do in Theorem \ref{C2leq5} below.
We begin by collecting (from Theorems \ref{thm:hyperelliptic},
\ref{thm:irreduciblen+reven}, and \ref{thm:irreduciblen+rodd})
the systems with minimal $C^2$, at most five.

\begin{proposition}\label{prop:minimalatmostfive}
Suppose that the full set of hypotheses hold.
Suppose that $C^2\leq 5$ and is the minimum possible for that $n$ and $r$.
Then the system is on the following list.
\[
\begin{array}{ccccccc}
\text{System} & \text{hyperelliptic?} & n & r & C^2 & g & g' \\
\hline
|9;3^8,2^2| & \text{hyp.} & 10 & 0 & 1 & 2  & 2 \\
|6;2^8,1^3| & \text{hyp.} & 11 & 0 & 1 & 2 & 1 \\
|4;2,1^{11}| & \text{hyp.} & 12 & 0 & 1 & 2 & 0 \\
\hline
%
|6;2^8,1^2| & \text{hyp.} & 10 & 1 & 2 & 2 & 1 \\
|4;2,1^{10}| & \text{hyp.} & 11 & 1 & 2 & 3 & 0 \\
|6;2^7,1^6| & \text{non-hyp.} & 13 & 0 & 2 & 3 & 1 \\
|4;1^{14}| & \text{non-hyp.} & 14 & 0 & 2 & 3 & 1 \\
|5;3,1^{14}| & \text{hyp.} & 15 & 0 & 2 & 3 & - \\
\hline
%
|4;2,1^9| & \text{hyp.} & 10 & 2 & 3 & 2 & 0 \\
|6;2^7,1^5| & \text{non-hyp.} & 12 & 1 & 3 & 1 & 1 \\
|4;1^{13}| & \text{non-hyp.} & 13 & 1 & 3 & 3 & 0 \\
|5;3,1^{13}| & \text{hyp.} & 14 & 1 & 3 & 3 & - \\
|5;2^2,1^{14}| & \text{non-hyp.} & 16 & 0 & 3 & 4 & 0 \\
|6;4,1^{17}| & \text{hyp.} & 18 & 0 & 3 & 4 & - \\
\hline
|6;2^7,1^4| & \text{non-hyp.} & 11 & 2 & 4 & 3 & 1 \\
|4;1^{12}| & \text{non-hyp.} & 12 & 2 & 4 & 3 & 0 \\
|5;3,1^{12}| & \text{hyp.} & 13 & 2 & 4 & 3 & - \\
|5;2^2,1^{13}| & \text{non-hyp.} & 15 & 1 & 4 & 4 & 0 \\
|6;2^5,1^{12}| & \text{non-hyp.} & 17 & 0 & 4 & 5 & 1 \\
|6;4,1^{16}| & \text{hyp.} & 17 & 1 & 4 & 4 & - \\
|7;5,1^{20}| & \text{hyp.} & 21 & 0 & 4 & 5 & - \\
\hline
|6;2^7,1^3| & \text{non-hyp.} & 10 & 3 & 5 & 3 & 1 \\
|4;1^{11}| & \text{non-hyp.} & 11 & 3 & 5 & 3 & 0 \\
|5;3,1^{11}| & \text{hyp.} & 12 & 3 & 5 & 3 & - \\
|5;2^2,1^{12}| & \text{non-hyp.} & 14 & 2 & 5 & 4 & 0 \\
|6;2^5,1^{11}| & \text{non-hyp.} & 16 & 1 & 5 & 1 & 1 \\
|6;4,1^{15}| & \text{hyp.} & 16 & 2 & 5 & 4 & - \\
|6;2^4,1^{15}| & \text{non-hyp.} & 19 & 0 & 5 & 6 & 1 \\
|5;1^{20}| & \text{non-hyp.} & 20 & 0 & 5 & 6 & 0 \\
|6;3,2,1^{18}| & \text{non-hyp.} & 20 & 0 & 5 & 6 & 0 \\
|7;5,1^{9}| & \text{hyp.} & 20 & 1 & 5 & 5 & - \\
|8;6,1^{23}| & \text{hyp.} & 24 & 0 & 5 & 6 & - \\
\end{array}
\]
(If $g'$ is indicated as '-', this means that the adjoint system is composed with a pencil.)

\end{proposition}

\begin{proof}
If the system is hyperelliptic, Theorem \ref{thm:hyperelliptic} gives a classification.
For the systems (i) there, we have $C^2 \geq 1$;
for the systems (ii), we have that $C^2=(n/3)+(4r/3)-3 \geq (10/3)-3 > 0$;
since the self-intersection of $|9;3^8,2^2|$ is $1$,
we have $C^2 \geq 1$ in all cases.

If $C^2=1$, we either have the system of degree $6$ with $n=11$ in case (i),
or in case (ii), we have $n+4r=12$,
which forces $n=12$ and $r=0$, hence $g=2$, 
giving $|4;2,1^{11}|$.
Case (iii) gives the system of degree $9$.

If $C^2=2$, the case (i) gives the system $|6;2^8,1^2|$;
the case (ii) forces $n+4r=15$, so we must have $(n,r)=(11,1)$ or $(15,0)$,
leading to the systems $|4;2,1^{10}|$ or $|5;3,1^{14}|$ indicated above.

If $C^2 \geq 3$, we must be in case (ii) of Theorem \ref{thm:hyperelliptic}.
If $C^2=3$ we then have $n+4r=18$, so $(n,r) = (10,2), (14,1), (18,0)$,
giving the three systems above. 
If $C^2=4$, we must then have $n+4r=21$, and the only possibilities are 
 $(n,r)=(13,2)$, $(17,1)$, $(21,0)$, giving the three systems in the list.

Finally if If $C^2=5$, then $n+4r=24$; the only possibilities are 
$(n,r)=(12,3)$, $(16,2)$, $(20,1)$, $(24,0)$
giving the four systems in the final section.

We will see (by listing the minimal self-intersection non-hyperelliptic systems below)
that all of these hyperelliptic systems give the minimal $C^2$ for their $(n,r)$
among all systems (hyperelliptic or not). 
Hence they all belong on the list above.

If the system is non-hyperelliptic with $n+r$ even,
then Theorem \ref{thm:irreduciblen+reven} and Lemma \ref{g'=0}
give the classification.
In case (a) of Lemma \ref{g'=0},
 the only system with minimum $C^2\leq 5$ is $|5;1^{20}|$,
which has $C^2=5$.

In case (b) of Lemma \ref{g'=0}
the $t=1$, $n=11, 12, 13, 14$ cases give the four systems
$|4;1^n|$ (with $C^2=16-n$ and $r=14-n$)
and the $t=2$, $n=17,18$ cases give the two systems $|5;2,1^{n-1}|$
(with $C^2 = 22-n$ and $r=18-n$);
all other cases have $C^2 \geq 6$.
For both of the $t=2$ cases however, 
there are hyperelliptic systems with smaller $C^2$,
so they do not go on the list.  The $t=1$ cases do.

In case (c) of the Lemma
the $t=2$, $n=14,15,16$ cases give the three systems
$|5;2^2,1^{n-2}|$ (with $C^2=19-n$ and $r=16-n$);
the $t=3$, $n=20$ case gives the system $|6,3,2,1^{18}|$
(with $C^2=5$ and $r=0$).
All other cases have $C^2\geq 6$.
All of these go on the list.
We note that the system of degree $6$ has the same parameters
($(n,r,C^2,g)=(20,0,5,6)$) as the system $|5;1^{20}|$ seen before;
here we have an example where the minimal system is not unique.

This produces all the systems of case (A) of Theorem \ref{thm:irreduciblen+reven};
case (B) adds no new non-hyperelliptic systems.

If the system is non-hyperelliptic with $n+r$ odd,
then Theorem \ref{thm:irreduciblen+rodd} gives the classification,
and we first analyze those systems coming from case (A)(a) of that theorem 
which have self-intersection $C^2 \leq 5$.

In that case (A)(a) we have $C^2 = 36-3\ell-n$ and since $n \leq 27-2\ell$,
we see that $C^2 \geq (36-3\ell)-(27-2\ell) = 9-\ell$.
This leads to ten of the systems given above, 
with $\ell \in\{4,5,6,7\}$, and $36-3\ell-n\leq 5$, 
or $31-3\ell \leq n \leq 27-2\ell$.
We see that $(\ell,n)$ is in the set 
$\{(4,19),(5,16),(5,17),
(6,13),(6,14),(6,15),
(7,10),(7,11),(7,12),(7,13)\}$.
The $(6,13)$, $(6,14)$, and $(6,15)$ cases lead to systems with $C^2=5,4,3$ respectively, and for these $(n,r)$ parameters 
we already have systems with $C^2=4,3,2$ respectively, 
of the form $|5;3,1^{n-1}|$.
The other seven cases lead to the seven non-hyperelliptic systems 
of degree $6$, with $n+r$ odd,
listed above that have only simple or double base points.

%
Finally we observe that in the other cases of Theorem \ref{thm:irreduciblen+rodd},
the minimal self-intersection is greater than five.	

In case (A)(b) of that Theorem the minimum $C^2$ is $8$.
In case (A)(c) the minimum is achieved with $\ell=4$, which forces $k \leq 23$;
the minimum $C^2$ is then $49-4\ell-k = 33-23 = 10$.
In case (A)(d) we get minimum $C^2$ by taking $\ell=1$ and $k = 33$, 
which gives $C^2 = 18$.	
In case (A)(e) the minimum $C^2$ is achieved with $m=4$; 
then $r \geq 0$ implies that $k\leq 10b+2$, 
so that the minimum $C^2$ for fixed $b$
is achieved with $k=10b+2$.  
This gives $C^2 =16b-16-10b-2+16 = 6b-2$, 
so we only have $C^2 \leq 5$ if $b=1$.
However the $b=1$ case is covered with the cases of (A)(a).

Finally the cases (B) of Theorem \ref{thm:irreduciblen+rodd} 
have either $(n,r)=(10,1)$ with $C^2=3$ or $(n,r)=(11,0)$ with $C^2=2$ 
so they do not have minimal selfintersection and don't belong in the list.

The genus of each system and of its adjoint are immediate to compute.
\end{proof}

We can now proceed to the classification of systems with $C^2 \leq 5$.
For the cases $C^2 \geq 3$, we will restrict ourselves to those systems with $r \geq 2$.
These systems are the ones that give dominant maps,
which we are primarily interested in.
The classification of the $r\leq 1$ systems is considerably more complicated.

Proposition \ref{prop:minimalatmostfive} gives us the systems with $C^2\leq 5$
that have minimal $C^2$ for their parameters $(n,r)$.
Our procedure for classifying the full set of systems with $C^2 \leq 5$
is to analyze systems with a given $C^2$ which are \emph{not} on the list above.
For such a system with a given $(n,r)$, if we are dealing with the case $C^2=s$,
then we can restrict ourselves to parameters $(n,r)$
which appear on the above list with self-intersection strictly less than $s$:
if the system is not on the list for that given $(n,r)$,
then it is not minimal for that $(n,r)$,
and hence the minimal system with that $(n,r)$ \emph{is} on the list
with a lower $C^2$.

\begin{theorem}\label{C2leq5}
Suppose the full set of hypotheses hold, 
in particular $n \geq 10$, $r \geq 0$, and $g \geq 2$.
Let $|C|$ be a complete linear system satisfying these, with $C^2 \leq 5$.
Then:
\begin{itemize}
\item[(a)] $C^2 \geq 1$, and if $C^2 =1$, 
then $|C|$ is hyperelliptic, 
equal to $|9;3^8,2^2|$, $|6;2^8,1^3|$ or $|4;2,1^{11}|$.
\item[(b)] If $C^2=2$, 
then either $|C|$ is on the Proposition \ref{prop:minimalatmostfive} list,
or is equal to one of the non-hyperelliptic systems 
$|15;5^7,4^3|$, $|18;6^8,5,3|$,
$|9;3^7,2^4|$, $|7;3,2^9,1^2|$,  or $|9;3^8,2,1^3|$.
\item[(c)] If $C^2=3$, and $r \geq 1$, 
then either $|C|$ is on the Proposition \ref{prop:minimalatmostfive} list,
or is equal to one of the systems $|10;4,3^9|$, $|12;4^8,3,2|$,
$|7;3,2^9,1|$, or $|9;3^8,2,1^2|$.
\item[(d)]
If $C^2 = 4$, and $r \geq 2$, 
then either $|C|$ is on the Proposition \ref{prop:minimalatmostfive} list,
or $C$ is equal to one of the systems
$|7;3, 2^9|$, or $|9;3^8,2,1|$.
\item [(e)] If $C^2=5$, and $r \geq 2$,
then either $|C|$ is on the Proposition \ref{prop:minimalatmostfive} list,
or $|C|$ is equal to one of the systems 
$|7;3, 2^8, 1^3|$, 
$|9;3^8,1^4|$, 
$|12;4^7,3^3|$,
$|15;5^8, 4, 2|$, 
$|6;2^6,1^7|$,
$|7;2^{11}|$, or
$|9;3^7,2^3,1|$.
\end{itemize}
These systems are collected in the list below.
\end{theorem}
\[
\begin{array}{ccccccc}
\text{System} & \text{hyperelliptic?} & n & r & C^2 & g & g' \\
\hline
|15;5^7,4^3| & \text{non-hyp.} & 10 & 0 & 2 & 3  & 4 \\
|18;6^8,5,3| & \text{non-hyp.} & 10 & 0 & 2 & 3  & 4 \\
|9;3^7,2^4| & \text{non-hyp.} & 11 & 0 & 2 & 3  & 3 \\
|7;3,2^9,1^2| & \text{non-hyp.} & 12 & 0 & 2 & 3  & 2 \\
|9;3^8,2,1^3| & \text{non-hyp.} & 12 & 0 & 2 & 3  & 2 \\
\hline
|10;4,3^9| & \text{non-hyp.} & 10 & 1 & 3 & 3  & 3 \\
|12;4^8,3,2| & \text{non-hyp.} & 10 & 1 & 3 & 3  & 3 \\
|7;3,2^9,1| & \text{non-hyp.} & 11 & 1 & 3 & 3  & 2 \\
|9;3^8,2,1^2| & \text{non-hyp.} & 11 & 1 & 3 & 3  & 2 \\
\hline
|7;3,2^9| & \text{non-hyp.} & 10 & 2 & 4 & 3  & 2 \\
|9;3^8,2,1| & \text{non-hyp.} & 10 & 2 & 4 & 3  & 2 \\
\hline
|12;4^7,3^3| & \text{non-hyp.} & 10 & 2 & 5 & 4  & 4 \\
|15;5^8,4,2| & \text{non-hyp.} & 10 & 2 & 5 & 4  & 4 \\
|7;2^{11}| & \text{non-hyp.} & 11 & 2 & 5 & 4  & 3 \\
|9;3^7,2^3,1| & \text{non-hyp.} & 11 & 2 & 5 & 4  & 3 \\
|7;3,2^8,1^3| & \text{non-hyp.} & 12 & 2 & 5 & 4  & 2 \\
|9;3^8,1^4| & \text{non-hyp.} & 12 & 2 & 5 & 4  & 2 \\
|6;2^6,1^7| & \text{non-hyp.} & 13 & 2 & 5 & 4  & 1 \\
\end{array}
\]

Statement (a) follows immediately by inspection of the classification given in Proposition \ref{prop:minimalatmostfive}.



We also note that if the system is hyperelliptic,
then the parameters $n$ and $r$ determine $C^2$;
hence the hyperelliptic systems of Proposition \ref{prop:minimalatmostfive}
are the only hyperelliptic systems with $C^2 \leq 5$.
Moreover, reviewing the proof above,
we conclude that every hyperelliptic system with $C^2 \leq 5$
is on that list;
it is never the case that a non-hyperelliptic system has a smaller self-intersection
for a given $n$ and $r$ than the hyperelliptic system (if it exists).

Therefore we may restrict the analysis to the non-hyperelliptic cases,
and thus also assume $g \geq 3$.


\begin{proof}[Proof (of Theorem \ref{C2leq5}(b)):]
If $|C|$ is non-hyperelliptic and not on the list,
then it is not minimal for its $(n,r)$;
hence the minimal system for its $(n,r)$ must have $C^2=1$,
and we conclude 
from the prior analysis of the minimal systems with $C^2=1$,
that $(n,r)=(10,0), (11,0)$ or $(12,0)$.

If $(n,r)=(10,0)$ then using Proposition \ref{prop:max} we have $g'=4$;
using property (F) of the Introduction we see that $g = 3$, so that $r'=2$, 
and $(C')^2=5$.
We claim that the system $|C'|$ is not hyperelliptic. 
Indeed, by Theorem \ref{thm:hyperelliptic} 
a hyperelliptic system with genus $4$ and dimension $0$ 
is Cremona equivalent to the $18$-point system $|6;4,1^{17}|$, 
whereas $|C'|$ only involves 10 points.
Therefore we can apply Lemma \ref{lem:jkl} below 
to deduce that this adjoint system
must be either $|12;4^7,3^3|$ or $|15;5^8,4,2|$.
Hence $|C| = |15;5^7,4^3|$ or $|18;6^8,5,3|$.

If $(n,r)=(11,0)$, then by Theorem \ref{thm:irreduciblen+rodd}(B),
the system must be $|9;3^7,2^4|$.

If $(n,r) = (12,0)$, then by Proposition \ref{prop:max} 
we have $n+3r+g' = 14$, so that $g'=2$.
In this case one has $g=3$ (by property (F))
and the adjoint system must have dimension $2$.
Using Proposition \ref{prop:classg=012}(2),
those giving systems of dimension $2$ are 
$|C'|=|4;2,1^9,0^{2}|$ and $|6;2^8,1,0^{3}|$.
These give the systems $|C|=|7;3,2^9,1^2|$ and $|9;3^8,2,1^3|$.
\end{proof}

\begin{proof}[Proof (of Theorem \ref{C2leq5}(c)):]
We note that if $C^2=3$ and is not on the list,
then $(n,r)$ must appear on the list for $C^2\leq 2$,
and since we are assuming $r \geq 1$, 
the only possibilities for $(n,r)$ are $(10,1)$ and $(11,1)$.

If $(n,r)=(10,1)$, we may invoke Theorem \ref{thm:irreduciblen+rodd}(B),
and conclude that $|C|$ is either $|10;4,3^9|$ or $|12;4^8,3,2|$.

If $(n,r)=(11,1)$, then $g=3$ and 
(using \eqref{eq:self} with $h=6$ and $\epsilon = 0$)
we have $g'=2$ and $r'=2$ for the nef adjoint system $|C'|$.
Moreover $(C')^2=3$.
Now we can use Proposition \ref{prop:classg=012}(2)
to conclude that $|C'|$ is either $|4;2,1^9,0|$ or $|6;2^8,1,0^2|$,
which then implies $|C|$ is either $|7;3,2^9,1|$ or $|9;3^8,2,1^2|$.

\end{proof}

\begin{proof}[Proof (of Theorem \ref{C2leq5}(d)):]


If $|C|$ is non-hyperelliptic and not on the list,
then it is not minimal for its $(n,r)$; hence the minimal system for its $(n,r)$
must have $C^2 \leq 3$ and $r \geq 2$.
We conclude from the list that the only possibility is $(n,r)=(10,2)$.
Hence $g=C^2+1-r=3$, and
then using \eqref{eq:self} we see that $h=6$ and $\epsilon=0$, so that $g'=2$.
We then can apply Proposition \ref{prop:classg=012}(2);
since $r' = g-1=2$,
we see that $|C'|$ is either $|4;2, 1^9|$ or $|6;2^8, 1,0|$.
Hence $|C|$ is either $|7;3,2^9|$ or $|9;3^8,2,1|$.

\end{proof} 

We will handle the $C^2=5$ cases via a series of lemmas.

\begin{lemma}\label{lem:n+revenc2=5}
If $C^2=5$, $r \geq 2$, and $n+r=2h$ is even,
then either $|C|$ is on the Proposition \ref{prop:minimalatmostfive} list,
or $|C|$ is equal to the system
$|7;3,2^8,1^3|$, $|9;3^8,1^4|$, 
or $n=10$ and $r=2$.
\end{lemma}

\begin{proof}
Again we must have $(n,r)$ on the list for lower self-intersection, and $r\geq 2$,
which implies that $(n,r)=(10,2)$ or $(12,2)$.
If $(n,r)=(12,2)$, then $g=4$ and \eqref{eq:self} gives $g'=2$ and $r'=3$.


Since $g'=2$ we can apply Proposition \ref{prop:classg=012}(2),
and deduce that $|C'|$ is either 
	$|4;2, 1^8,0^2|$ or  
	$|6;2^8, 0^4|$;
these lead to $|C|$ being either $|7;3,2^8,1^3|$ or $|9;3^8,1^4|$.
%
%
%
%
\end{proof}

The case $n=10$, $r=2$ in the above Lemma requires some special treatment,
to which we now turn.

\begin{lemma}\label{lem:jkl} Suppose that the system $|C|$ is non-hyperelliptic 
	with $C^2 = 5$, $n=10$, $r=2$.
	Then $C$ is either Cremona equivalent to $|12;4^7,3^3|$ or to $|15;5^8,4,2|$.  
\end{lemma}

\begin{proof} 
Since $C^2=5$ and $r=2$, then $g=C^2-r+1=4$ 
by observation \eqref{hypothesis:cohomology}.	
Moreover we have $C\cdot K_{10}=1$. 
	
	As we saw in the proof of Lemma \ref {lem:n+revenc2=5}, 
	we have also $g'=4$, and the general curve in $|C'|$ is irreducible
		by the assumption that $|C|$ is non-hyperellyptic.
	One has $(C')^2=(C+K_{10})^2=6$.
	
	Consider the second adjoint system $|C''|=|C+2K_{10}|$, 
	and its Zariski decomposition 
	$$
	C+2K_n\sim P+ A
	$$
	in which, as seen in Lemma \ref {lem:a}, $A=A_1+\dots+A_h$
	where each $A_i$ is a $(-1)$-curve with $C\cdot A_i=1$.
	One has $\dim (|P|)=g'-1=3$. 
	
	\begin{claim} \label{cl:1}  
		The general curve in $|P|$ is irreducible. 
	\end{claim}
	
	\begin{proof}[Proof of Claim \ref {cl:1}] 
		Suppose the claim is not true, i.e.,  
		$|P|$ is composed with a pencil $|F|$ whose general curve is irreducible. 
		Then since $\dim|P|=3$, we have $P\sim 3F$ and $C'\cdot F=2$ 
		and the general curve in $|C'|$ is hyperelliptic. 
		Moreover, since $(C')^2=6$, by the index theorem one has $F^2=0$.

		Consider the exact sequence
		$$
		0\longrightarrow \mathcal O_{X_{10}}(C'-F)\longrightarrow \mathcal O_{X_{10}}(C')\longrightarrow \mathcal O_{F}(C')\longrightarrow 0.
		$$
		We have $h^0(X_{10}, \mathcal O_{X_{10}}(C'))=g=4$ and 
		$h^0(F, \mathcal O_{F}(C'))\leq 3$, because $C'\cdot F=2$. 
		Hence $h^0(X_{10}, \mathcal O_{X_{10}}(C'-F))>0$, i.e., $C'-F$ is effective 
		and therefore $h^2(X_{10}, \mathcal O_{X_{10}}(C'-F))=0$. 
		Moreover $h^1(X_{10}, \mathcal O_{X_{10}}(C'))=0$. 
		This implies that $h^1(F, \mathcal O_{F}(C'))=0$. 
		On the other hand $h^0(F, \mathcal O_{F}(C'))\geq 2$, 
		since otherwise $|C'|$ would have some fixed point on the general curve of $|F|$ 
		and therefore $|C'|$ would have some fixed component, a contradiction. 
		So, if $f$ is the genus of the curves in $|F|$, 
		by the Riermann--Roch theorem we have 
		$2\leq h^0(F, \mathcal O_{F}(C'))=3-f$, i.e., $f\leq 1$. 
		
		If $f=0$, then $|F|$ is a pencil of rational curves, 
		and therefore it is Cremona equivalent to $|1;1|$. 
		Then $C''$ is Cremona equivalent to $|3;3,0^{9-h}, (-1)^h|$ 
		and $|C|$ is Cremona equivalent to $|9;5, 2^{9-h}, 1^h|$,
		contradicting $C^2=5$.
		
		If $f=1$, then $|F|$ is Cremona equivalent to $|3;1^8|$,  
		$C''$ is Cremona equivalent to $|9;3^8,0^{2-h}, (-1)^h|$ 
		and $|C|$ is Cremona equivalent to $|15;5^8, 2^{2-h}, 1^h|$ 
		again contradicting $C^2=5$. 
	\end{proof}
	
	\begin{claim} \label{cl:2} 
		One has $h=0$, i.e., $|C''|$ has no fixed components 
		and therefore its general curve is irreducible. 
	\end{claim}
	
	\begin{proof} [Proof of Claim \ref {cl:2}] 
		Suppose by contradiction that $h\geq 1$. 
		Then $C'\cdot A_i=0$, for $1\leq i\leq h$. 
		We can then contract the $(-1)$--curves $A_1,\ldots, A_h$, 
		so that we can work on $X_{10-h}$, 
		and we  denote by $\bar C'$ the image of $C'$ on $X_{10-h}$. 
		We  have $(\bar C')^2=(C')^2=6$ and the genus of $\bar C'$ is $4$, 
		so that $\bar C'\cdot K_{10-h}=0$. 
		By the index theorem this implies that $K_{10-h}^2\leq 0$ 
		which implies $h=1$ and therefore $h^0(X_9, \mathcal O_{X_9}(-K_9))=1$, 
		i.e., there is a unique anticanonical curve $D$ on $X_9$. 
			Now, it is well known that $D$ is the only irreducible effective 
			curve on $X_9$ with $D\cdot K_9=0$, so $\bar C'=D$, a contradiction
			with $g(\bar C')=4$.
	\end{proof}
	
	Since $C''\cdot (C''+K_{10})=(C+2K_{10})\cdot (C+3K_{10})=4$, 
	we see that the curves in $|C''|$ have (arithmetic) genus $g''=3$. 
	
	Next we consider the third adjoint system $|C'''|=|C+3K_{10}|$ 
	which has dimension $g''-1=2$,	
		and its Zariski decomposition 
		$$
		C+3K_n\sim P'+ A'.
		$$
		with $A'=A_1+\dots+A_h$
		where each $A_i$ is a $(-1)$-curve with $C\cdot A_i=-2$.

	\begin{claim} \label{cl:3} The general curve in $|P'|$ is irreducible.
	\end{claim}
	
	\begin{proof} [Proof of Claim \ref {cl:3}] 
		The proof is similar to the one of Claim \ref {cl:1}, so we will be brief. 
		
		Suppose by contradiction that $|P'|$ is composed with a pencil $|G|$ 
		whose general curve is irreducible. 
		Then the general curve in $C''$ is hyperelliptic, 
		so that $P'\sim 2G$ and $G\cdot C''=2$. 
		Moreover, since $(C'')^2=(C+2K_{10})^2=5$, by the index theorem one has $G^2=0$.
		
		Consider the exact sequence
		$$
		0\longrightarrow \mathcal O_{X_{10}}(C''-G)\longrightarrow \mathcal O_{X_{10}}(C'')\longrightarrow \mathcal O_{G}(C'')\longrightarrow 0.
		$$
		We have $h^0(X_{10}, \mathcal O_{X_{10}}(C''))=4$ and 
		$h^0(G, \mathcal O_{G}(C''))\leq 3$, because $C''\cdot G=2$. 
		Hence $h^0(X_{10}, \mathcal O_{X_{10}}(C''-G))>0$,  
		thus $h^2(X_{10}, \mathcal O_{X_{10}}(C''-G))=0$. 
		Moreover $h^1(X_{10}, \mathcal O_{X_{10}}(C''))=0$. 
		This implies that $h^1(G, \mathcal O_{G}(C''))=0$. 
		On the other hand $h^0(G, \mathcal O_{G}(C''))\geq 2$. 
		So, if $\gamma$ is the genus of the curves in $|G|$, 
		by the Riemann--Roch theorem we have $2\leq h^0(G, \mathcal O_{G}(C''))=3-f$, 
		i.e., $\gamma\leq 1$. 
		
		If $\gamma=0$, then $|G|$ is  Cremona equivalent to $|1;1|$. 
		Then $|C'''|$ is Cremona equivalent to $|2;2,0^{9-h}, (-1)^h|$ 
		and then $|C|$ is Cremona equivalent to $|11;5, 3^{9-h}, 2^h|$
		which contradicts $C^2=5$.
		
		If $\gamma=1$, then $|G|$ is Cremona equivalent to $|3;1^8|$. 
		Then $|C'''|$ is Cremona equivalent to $|6;2^8,0^{2-h}, (-1)^h|$ 
		and then $|C|$ is Cremona equivalent to $|15;5^8, 3^{2-h}, 2^h|$
		and again this contradicts $C^2=5$. 
	\end{proof}
	
	\begin{claim} \label{cl:4} If $|C'''|$ has some fixed $(-1)$--curve, then it has only one such fixed curve. \end{claim}
	
	\begin{proof} [Proof of Claim \ref {cl:4}] 
		Suppose  that $|C'''|$ has some fixed $(-1)$--curve, 
			i.e., 	 $A=A_1+\dots+A_h$ with $h>0$.
		Then $C''\cdot A_i=0$, for $1\leq i\leq h$. 
		We can  contract the $(-1)$--curves $A_1,\ldots, A_h$, 
		so that we can work on $X_{10-h}$, 
		and we  denote by $\bar C''$ the image of $C''$ on $X_{10-h}$. 
		We  have $(\bar C'')^2=(C'')^2=5$ and the genus of $\bar C''$ is $3$, 
		so that $\bar C''\cdot K_{10-h}=-1$. 
		If $h>2$, then $\dim (|{-K_{10-h}}|)\geq 2$ and $|{-K_{10-h}}|$ is base point free. 
		Since $\bar C''\cdot ({-K_{10-h}})=1$ and $g''=3$, this is impossible.
		Hence $h\leq 2$. 
		
		If $h=2$, one has $\dim (|{-K_{8}}|)=1$ and $|{-K_{10-h}}|$ has a base point $x$. 
		Then, by the same argument as above, $x$ has to be also a base point for $|\bar C''|$. 
		If we blow--up $x$, we get a surface $X$, 
		with the strict transforms $|\tilde C''|$ of $|\bar C''|$ and $|D|$ of the pencil  $|{-K_{8}}|$.  
		One has $(\tilde C'')^2=4$, $D^2=0$ and $\tilde C''\cdot D=0$, 
		which violates the index theorem. 	
	\end{proof}
	
	So we must have $h\leq 1$.
	
	Since $C'''\cdot (C'''+K_{10})=(C+3K_{10})\cdot (C+4K_{10})=0$, 
	the curves in $|C'''|$ have arithmetic genus $1$. 
%
%
%
		But $P'\cdot A_i=0$ in the Zariski decomposition, so
	if $h=1$ we have $1=p_a(A_1)+p_a(P')-1=p_a(P')-1$, 
	hence $p_a(P')=2$. 
	So we can apply Proposition \ref{prop:classg=012}(2)
to the linear system  $|P'|$ of dimension 2,  
	and we deduce that it is Cremona equivalent to of one of the following:
	$|4;2, 1^9|$,  $|6;2^8, 1, 0|$. 
	However $|4;2, 1^9|$ is not possible 
	because this system intersects positively any $(-1)$--curve. 
	Indeed, if not, it would be Cremona equivalent to a linear system on $X_n$ with $n\leq 9$, which is not possible by \cite [Prop. 10.10]{CCi}. 
	So we can assume that $|P'|$ is $|6;2^8, 1, 0|$ and $|C'''|$ is $|6;2^8, 1, -1|$. 
	Accordingly $|C|$ is then Cremona equivalent to $|15;5^8,4,2|$. 
	
	Finally, we can assume that $h=0$. 
	Then the general curve in $|C'''|$ is irreducible of arithmetic genus 1, 
	so that $|C'''|$ is Cremona equivalent to $|3;1^7,0^3|$, 
	so that $|C|$ is Cremona equivalent to $|12;4^7,3^3|$.   
	
This finishes the proof of Lemma \ref{lem:jkl}.
\end{proof}

%

This completes the $C^2=5$, $n+r$ even case;
we now turn our attention to the $n+r$ odd case.

\begin{lemma}\label{lem:n+rodd}
If $C^2=5$, $r \geq 2$, and $n+r=2h+1$ is odd,
then either $|C|$ is on the Proposition \ref{prop:minimalatmostfive} list,
or $|C|$ is equal to one of the two-dimensional systems
$|6;2^6,1^7|$,
$|7;2^{11}|$, or
$|9;3^7,2^3,1|$.
\end{lemma}

\begin{proof}
%
Since we are assuming $|C|$ is not on the list of Proposition \ref{prop:minimalatmostfive},
then the $(n,r)$ values must occur with the smaller $C^2$ cases ($\leq 4$);
this means that $(n,r) = (11,2)$ or $(13,2)$ since $r\geq 2$.
	If $C^2=5$, by Proposition \ref {prop:max} we must have $n+3r+g'=20$
	so $g=14-n$, and the possible cases are $(n,r,g') = (11,2,3), (13,2,1)$.

	If $(n,r,g') = (13,2,1)$, we apply Lemma \ref {g'=1}. 
	Case (a) of that lemma gives us the linear system 
$|6;2^6,1^7|$. 
	In case (b) one has no contribution. 
	
	In the case	$(n,r,g') = (11,2,3)$, 
	consider the system $|C''|=|2K_n+C|$ to $|C'|$, 
	that has dimension $g'-1=2$. 
If $|C''|$ is nef, then $|C''|$ has no fixed component (see Proposition \ref {prop:g2}). If $|C''|$ is not nef, consider its Zariski decomposition 
		$$
		C+2K_n\sim P+ A
		$$
in which, as it follows by Lemma \ref {lem:a}, $A=A_1+\dots+A_h$
where each $A_i$ is a $(-1)$-curve with $C\cdot A_i=1$ 
and $A_i\cdot A_j=0$ if $1\leq i<j\leq h$.
	So $|C''|$ is non--special, 
	and it may have some $(-1)$--curve in its fixed part. 
	By blowing down these $(-1)$--curves, 
	we may assume that  $|C''|$ has no fixed part and it lives on $X_m$, 
	with $m\leq n$ (see \cite [Lemma 5.1]{CM}). 
	
	Suppose first the general curve in $|C''|$ is reducible. 
	Then, by Proposition \ref{prop:fin}, 
	$|C''|$ is Cremona equivalent to $|2;2|$, 
	and therefore $|C'|$ is Cremona equivalent to $|5;3, 1^s|$ 
	and $|C|$ is Cremona equivalent to $|8;4, 2^s, 1^t|$,  with $s+t=10$. 
	Since $C^2=5$, we must have $43=4s+t$ which, together with $s+t=10$ 
	gives $s=11$, a contradiction. 
	So the general curve in $|C''|$ is irreducible. 
	Since $C'\cdot C''=4$, from the exact sequence
	$$
	0\longrightarrow \mathcal O_{X_n}(-K_n)\longrightarrow
	\mathcal O_{X_n}(C')\longrightarrow \mathcal O_{C''}(C')\longrightarrow 0
	$$ 
	we deduce that $|C'|$ (that has dimension $g-1=3$) 
	cuts out a $g^3_4$ on the general curve $C''$ of $|C''|$, 
	and therefore $C''$ has either genus 0 or genus 1. 
	
	If $C''$ has genus 0, then $|C''|$ is Cremona equivalent to the linear system $|1|$, 
	so that $|C'|$ is Cremona equivalent to $|4;1^s|$ 
	and $|C|$ to $|7; 2^s,1^t|$ with $s+t=11$. 
	Since $49-4s-t=C^2=5$, we find $s=11, b=0$, 
	hence we get the linear system $|7;2^{11}|$. 
	
	If $C''$ has genus 1, then $|C''|$ is Cremona equivalent to $|3;1^7|$, 
	so that $|C'|$ is Cremona equivalent to $|6;2^7, 1^s|$ 
	and $|C|$ to $|9; 3^7,2^s, 1^t|$ with $s+t=4$. 
	Since $18-4s-t=C^2=5$, we find $s=3, t=1$, 
	hence we get the linear system $|9;3^7,2^{3}, 1|$. 
\end{proof}

\begin{proof}[Proof (of Theorem \ref{C2leq5}(e)):]
Statement (e) now follows from Lemma \ref{lem:n+revenc2=5},
Lemma \ref{lem:jkl}, and Lemma \ref{lem:n+rodd}.
\end{proof}

This completes the proof of Theorem \ref{C2leq5}.

\subsection{Minimal degree maps to the plane}
In the previous subsection we classified systems with low $C^2$ (given $(n,r)$).
We are also in a position to classify systems of minimal degree with low $r$.
Since we are mainly interested in dominant maps,
we present the results below for systems with $r=2$.

\begin{proposition}
Suppose that the full set of hypotheses are true.
Then the complete linear systems with $n\geq 10$ and $r=2$
of minimal self-intersection are the following.
\begin{itemize}
\item[(a)] If the adjoint system $|C'|$ is composed with a pencil,
	then $n \equiv 1$  (mod. 3) and $|C|$ is Cremona equivalent to 
	$|(n+2)/3; (n-4)/3, 1^{n-1}|$.
\item[(b)] If the adjoint system has  irreducible general member
	and $n$ is even, then $n \geq 12$ and the minimum self-intersection is $n/2-2$.
	Any such system achieving this minimum is Cremona equivalent to either
	\begin{enumerate}
	\item $|5;1^{18}|$ (here $n=18$),
	\item $|n/4+1;n/4-2,1^{n-1}|$ (here $4|n$, $n \geq 12$),
		or 
	\item $|(n+6)/4; (n-6)/4,2,1^{n-2}|$ (here $n \equiv 2 (mod. 4)$, $n \geq 10$).
	\end{enumerate}
\item[(c)] If the adjoint system has an irreducible general member and $n$ is odd,
	then the system with minimum self-intersection is Cremona equivalent to either:
	\begin{enumerate}
	\item $|6;2^\ell,1^{25-3\ell}|$ (here $\ell \in \{0,\ldots,7\}$, $n=25-2\ell$),
	\item $|7;3^2,1^{21}|$ (here $n=23$),
	\item $|7;2^\ell,1^{33-3\ell}|$ (here $\ell \in \{0,\ldots,4\}$, $n=33-2\ell$),
	\item $|8;3, 2^\ell,1^{36-3\ell}|$ (here $\ell \in\{0,1\}$, $n=37-2\ell$),
			or
	\item $|2b+4;2b,2^\ell,1^{10b-3\ell+12}|$ (here $\ell \in \{0,\ldots,4\}$).
	\end{enumerate} 
\end{itemize}
\end{proposition}

We note that the system $|4;2,1^9|$ with $r=2$ and $n=10$
occurs as the smallest system in (a),
and the smallest system in (b)(3).
It also occurs with  $b=0$ case in (c)(5), with $\ell=1$.

\begin{proof}
Statement (a) are the systems of Theorem \ref{thm:hyperelliptic}(ii),
which are the only ones there that have $r=2$ and $n\geq 10$.
Statement (b) are the systems of Lemma \ref{g'=0}, with $r=2$,
coming from (A) of Theorem \ref{thm:irreduciblen+reven}.
The list of statement (c) follows from (A) of Theorem \ref{thm:irreduciblen+rodd}.
\end{proof}

\subsection{Gaps} 

In this section we consider the problem of classifying the linear systems $|C|$ on $X_n$, 
with $r=\dim(|C|)\geq 2$ and with $C^2$ one more than the minimum. 
In this case the adjoint linear system $|C'|$ has irreducible general member,
and we are in the non-hyperelliptic case. 
As we will see, again  there will here be a substantial difference 
between the cases $n+r$ even or odd.

\subsubsection{The $n+r$ even case} 
First we examine the case $n+r$ even.  
Recall that we are assuming $g\geq 2$ (recall hypothesis \eqref{hypothesis:g}), 
so that the dimension of the adjoint system $|C'|$ is $g-1\geq 1$. 

\begin{proposition}\label{prop:cclass}
We assume that we have the full set of hypotheses.
In addition, suppose $|C|$ is a non-hyperelliptic system, with $r \geq 2$,
and $n+r=2h$ is even.
If  $|C|$  has self--intersection one more than the minimum, i.e.,  $C^2=h+r-4$, 
then $|C|$ is Cremona equivalent to one of the following:
\begin{itemize}
	\item[(a)] $|7;3, 2^k, 1^\ell|$ 
with $r=29-3k-\ell \geq 2$, $h=15-k$ and $C^2=40-4k-\ell$; 
	\item [b)] $|9;3^8, 2^k,1^ell|$ 
with $r=6-3k-\ell \geq 2$, $h=7-k$ and $C^2=9-4k-ell$.
	\end{itemize}
\end{proposition}

\begin{proof} 
From \eqref {eq:self} we see that $C^2$ is one more than the minimum 
if and only if $g'=2$. 
Then we apply Proposition \ref{prop:classg=012}(2),  
and we deduce that $|C'|$ can be assumed to be of one of the following types:
$|4;2, 1^k,0^{n-k-1}|$ with $k\leq 10$ (leading to case (a); or
$|6;2^8, 1^k,0^{n-k-8}|$ with $k\leq 2$ (leading to case (b). 
\end{proof}

\begin{remark}\label{rem:lki} 
	In all cases listed in  Proposition \ref {prop:cclass}, 
	one has that $n+r$ is bounded above. 
	This shows that for $n+r$ large enough (actually $n+r\geq 32$ is enough) 
	there is no linear system $|C|$ of dimension $r$ and genus $g\geq 2$ on $X_n$, 
	such that the adjoint linear system $|C'|$ has irreducible general member, 
	with self--intersection one more than the minimum.  
	This indicates that, for fixed $n$, $r$, 
	one may expect gaps in the possible values of $C^2$. 
	It would be interesting to determine all such gaps, 
	or at least the gaps for a given $n$ and $r\gg 0$. 
\end{remark}

\subsubsection{The $n+r$ odd case} We  focus now on  the case $n+r$ odd.  
We keep the notation introduced 
in the proof of Proposition \ref {prop:CEboundC'} and in Remark \ref {rem:for}.

\begin{proposition}\label{prop:lop}

We assume that we have the full set of hypotheses.
In addition, suppose $|C|$ is a non-hyperelliptic system, with $r \geq 2$,
and $n+r=2h+1$ is odd.
%
If  $|C|$  has self--intersection one more than the minimum, i.e.,  $C^2=h+b+r-4$,
(so that we must have $g'\geq 3$), 
then $|C|$ is Cremona equivalent to one of the following:
%
\begin{itemize}
\item [(i)] $|7;2^s, 1^t|$, with $5\leq s\leq 9$ and $r=35-3s-t\geq 2$;
\item [(ii)] $|9;3^s, 2^t, 1^u|$, with $0\leq s\leq 7$ odd, $t\leq 1$, and 
		$r=54-6s-3t-u\geq 2$;
\item [(iii)] $|11;5^2, 1^t|$ with $r=47-t\geq 2$,
		or $|10;4^2, 1^s|$ with $r=45-s\geq 2$;
\item [(iv)] $|g'+5;g'+1,2^s, 1^t|$, with $5\leq s\leq 9$ and $r=5g'+19-3s-t\geq 2$; 
\item [(v)] $|b+6; b+1, 2^s, 1^t|$ with $1\leq b\leq 4$, $6-2b\leq s\leq 10-2b$ and $r=6b+26-3s-t\geq 2$. 
\end{itemize}
All values of $(n,r)$ with $h\geq 8$ as above occur on this list.
\end{proposition}

\begin{proof} 
	Taking into account Proposition \ref {prop:CEboundC'} and  Remark \ref {rem:for}, 
	we have that $C^2$ is one more than the minimum if and only if $a=b+1$, 
	so that $C^2=h+b+r-4$, $g'=2b+1$, $g=h+b-3$. 
	Then we have
	\[
	4g'+4-(C')^2=3g'+5-\dim (|C'|)=5b+12-h=m+5.
	\]
	Since $0\leq m\leq 4$, we have 
	\[
	4g'-5\leq (C')^2\leq 4g'-1.
	\]
	As we saw in the proof of Theorem \ref{thm:irreduciblen+rodd}, 
	the map determined by the adjoint system $|C+K_n|$ 
	is a birational morphism of $X_n$ to its image. 
	So we can apply \cite[Thm. 1.2 and \S 3.1] {CD} 
	and the results in \cite [\S\S 9--10] {CDM}. 
	Precisely $|C'|$ must be Cremona equivalent to one of the following systems:
	\begin{itemize}
		\item [(i)] $|4;1^s|$, with $5\leq s\leq 9$;
		\item [(ii)] $|6;2^s, 1^t|$, with $s\leq 7$ odd and $t\leq 1$;
		\item [(iii)] $|8;4^2|$ or $|7;3^2|$ ;
		\item [(iv)] $|g'+2;g',1^s|$, with $5\leq s\leq 9$, or 
		$|g'+3;g'+1,2,1^s|$, with $5\leq s\leq 9$ and $g'$ odd
		(although the degree $g'+3$ system leads to a system $|C|$ which is not Cremona minimal); 
		\item [(v)] $|b+3; b, 1^s|$ with $1\leq b\leq 4$ and $6-2b\leq	 s\leq 10-2b$.
	\end{itemize}
	
	Accordingly, $|C|$ must be Cremona equivalent 
	to one of the systems listed in the statement of the proposition. 
	
Finally we have to prove the final statement that all $(n,r)$ are possible. 
Given $n$ and $r$, we have $n+r=2h+1$ and assume $h\geq 8$. 
Set $h-7=5b-m$, with $b\geq 1$ and $0\leq m\leq 4$. 
Fix $s=m+5\leq 9$ and $t=n-s-1$, so that $h=5b+12-s$. 
Then 
	\[
	n+r=2h+1=10b+25-2s, \quad \text{hence}\quad r=10b+25-2s-n.
	\]
Set now $g'=2b+1$ and consider the linear system $|g'+5;g'+1,2^s, 1^t|$ as in (iv). 
This system lives on $X_n$ with $n=s+t+1$ and $r=5g'+19-3s-t=10b+25-2s-n$ as desired. \end{proof}

\begin{remark}\label{rem:sgo} 
	Proposition \ref {prop:lop} does not apply to the cases $h\leq 7$, i.e., $n+r\leq 15$,
	which are the cases $(10+t,3-t)$, with $0\leq t\leq 1$, 
	and  $(10+s,5-s)$, with $0\leq s\leq 3$. 
	For these pairs $(n,r)$ one more than the minimum self--intersection 
is in fact a gap in the possible self-intersections. 
	
	On the other hand Proposition \ref {prop:lop} shows that 
	if $n+r=2h+1$ with $h\geq 8$, contrary to the case $n + r$ even, 
	the value of the self--intersectiom one more than the minimum is never a gap.
\end{remark}

\subsection{Base point freeness, birationality,  etc.}

In this section we will look at the linear systems on $X_n$, 
with $n\geq 10$ of dimension $r\geq 2$, 
with minimal self--intersection,  with regard to
base point freeness, birationality, etc.
We again assume the full set of hypotheses.
This implies that $C^2>0$. 
If $E$ is a negative curve, 
we have $C\cdot E>0$ by property (C). 
If $E$ has $E^2\geq 0$ and $C\cdot E=0$ we have a contradiction to the Index Theorem.
Hence by the Nakai--Moishezon Criterion all systems $|C|$ are ample
with our hypotheses.

To study base point freeness and birationality we need some preliminaries. 

Let $S\subset \mathbb P^r$, with $r\geq 3$, 
be an irreducible, projective, non--degenerate surface. 
If $p\in S$ is a general point, 
we can consider the projection $\pi: S\dasharrow S'\subset \mathbb P^{r-1}$ 
with center the point $p$, 
with $S'$ the image of this map. 
This is called a \emph{general internal projection} of $S$.  

If $\mathcal H$ is the linear system cut out on $S$ 
by the hyperplanes of $\mathbb P^r$, 
then  the  general internal projection 
$\pi: S\dasharrow S'\subset \mathbb P^{r-1}$ 
is the map determined by the linear system $\tilde {\mathcal H}$ 
on the blow--up $\tilde S$ of $S$ at $p$, 
which is the strict transform on $\tilde S$ of the system $\mathcal H(-p)$ 
of the linear system of hyperplanes sections of $S$ containing $p$. 

\begin{lemma}\label{lem:dop} 
	In the above set--up, 
	we have that $\tilde {\mathcal H}$ is base point free and, 
	if $r\geq 4$, 
	the general internal projection $\pi: S\dasharrow S'\subset \mathbb P^{r-1}$ 
	of $S$ from $p$ 
	is birational onto its image.
\end{lemma}

\begin{proof} The first assertion is trivial, 
	since the scheme-theoretical  intersection 
of all the hyperplanes in $\bbP^r$ 
	containing $p$ is $p$ itself. 
	As for the second assertion it follows from the well known fact that, if $r\geq 4$, 
	then the general secant line to $S$ is not trisecant 
	(this is the so--called \emph{Trisecant Lemma}, see \cite {L}). 
\end{proof}

\begin{proposition}\label{prop:bpf}
Assuming the full set of hypotheses, 
the linear systems on $X_n$, with $n\geq 10$, of dimension $r\geq 2$, 
with minimal self--intersection in the hyperelliptic and in the non--hyperelliptic case, 
are base point free and, if $r\geq 3$, 
they determine birational morphisms onto their images. 
\end{proposition}

\begin{proof} 
	The proof consists in a case by case analysis. 
	Let us start with the hyperelliptic case examined in Theorem \ref {thm:hyperelliptic};
	we need only consider case (ii) there, since $n \geq 10$ and $r \geq 2$. 
	Consider first the linear system $|g+2;g|$, with $g\geq 2$, on $X_1$. 
	This is clearly very ample, 
	and the map determined by it embeds $X_1$ in $\mathbb P^{3g+5}$ 
	as a smooth surface of degree $4g+4$. 
	A linear system of the form $|g+2;g,1^{n-1}|$ 
	corresponds to $n-1$ general internal projections of the above surface. 
	Hence the assertion follows right away from Lemma \ref {lem:dop}.
	
	Consider next the non--hyperelliptic case with $n+r$ even, 
	examined in Theorem \ref {thm:irreduciblen+reven}(A)
	and Lemma \ref {g'=0}. 
	The linear system $|5|$ is very ample on $\mathbb P^2$, 
	and the map determined by it embeds $\mathbb P^2$ in $\mathbb P^{20}$ 
	as a smooth surface of degree $25$. 
	The linear systems $|5;1^n|$ in (a) of Lemma \ref {g'=0}
	correspond to $n$ general internal projections of the above surface 
	and again the assertion follows  from Lemma \ref {lem:dop}. 
	The linear system $|t+3;t|$, $t\geq 1$, on $X_1$ is very ample, 
	and the map determined by it embeds $X_1$ in $\mathbb P^{3t+6}$ 
	as a smooth surface of degree $6t+9$. 
	A linear system of the form $|t+3;t, 1^{n-1}|$, 
	as in (b) of Lemma \ref {g'=0}, 
	corresponds to $n-1$ general internal projections of the above surface 
	and we get the same conclusion as above. 
	The linear system $|t+3;t,2|$, $t\geq 1$, on $X_2$ is easily seen to be very ample, 
	and the map determined by it embeds $X_2$ in $\mathbb P^{3t+3}$ 
	as a smooth surface of degree $6t+5$. 
	A linear system of the form $|t+3;t, 2,1^{n-2}|$, 
	as in (c) of Lemma \ref {g'=0}, 
	corresponds to $n-2$ general internal projections of the above surface 
	and we conclude  as above.
	
	Finally we focus on the non--hyperelliptic case with $n+r$ odd 
	considered in Theorem \ref {thm:irreduciblen+rodd}(A). 
	The approach is similar to the above and therefore we will be brief. 
	We simply have to check that the following linear systems 
	are very ample on the corresponding surfaces: 
	$|6;2^\ell|$ with $\ell \in \{0,\ldots,7\}$; 
	$|7;3^2|$; 
	$|7;2^\ell|$ with $0\leq \ell \leq 4$; 
	$|8;3,2^\ell|$ with $0\leq \ell\leq 1$; 
	$|2b+4;2b,2^m|$ with $0\leq m \leq 4$. 
	Each of these systems is clearly very ample.
\end{proof}

\begin{example}\label{ex:10} To  illustrate  the ideas, 
	let us consider the case $n=10$ and 
	look at the system of minimal self--intersection for low values of $r$, 
	i.e., $3\leq r\leq 5$. 
	
If $r=3$, we have the system $|6;2^7,1^3|$ with $g=3$ and self--intersection $5$. 
It determines a birational morphism of $X_{10}$ 
to a degree $5$ surface $S$ in $\mathbb P^3$. 
It has been proved in \cite [pp. 487--88]{Con} 
that the surface $S$ has double curve consisting of three distinct lines 
intersecting in one point that is triple for $S$. 
We can directly check this in the following way,
which provides a converse to the statement.

Consider a general quintic surface $S$ in $\mathbb P^3$ 
with a double curve consisting of three distinct lines 
intersecting in one point that is triple for $S$. 
Let $X\longrightarrow S$ be the normalization, with $X$ smooth. 
Let $G$ be the total transform on $X$ of the double curve. 
Let $C$ be the proper transform on $X$ of a general plane section of $S$. 
Note that $K_X\sim C-G$, 
so that $p_g(X)=0$ and clearly also $P_2(X)=0$. 
We claim that $X$ is regular 
and therefore it is rational by Castelnuovo Criterion of rationality. 
Indeed, let $|C'|$ be the adjoint system of $|C|$. 
Look at the exact sequence
	\[
	0\longrightarrow \mathcal O_X(K_X)\longrightarrow \mathcal O_X(C')\longrightarrow \mathcal O_C(C')\cong \omega_C\longrightarrow 0.
	\] 
Since $C'\sim 2C-G$, then $|C'|$ is the pull--back on $X$ 
of the linear system cut out on $S$ 
by the quadric cones that contain the double curve.
Then it is immediate that the map 
	$H^0(X, \mathcal O_X(C'))\longrightarrow H^0(C,\omega_C)$ is surjective, 
	and that $h^1(X, \mathcal O_X(C'))=0$ by the Kawamata--Viehweg theorem. 
	This implies $h^1(X, \mathcal O_X(K_X))=0$, hence $q(X)=0$. 
	
The description of $|C'|$ makes it easy to check that 
the genus of the curves in $|C'|$ is 1,
as is indicated in the list of Proposition \ref{prop:minimalatmostfive} 
and so we get $0=(2K_X+C)\cdot (K_X+C)$. 
Since $K_X\cdot C=-1$ and $C^2=5$ 
we deduce that $K_X^2=-1$, so that $X=X_{10}$. 
Then the normalization morphism $X_{10}=X\longrightarrow S$ 	
is determined by a linear system of dimension $3$ 
of curves with self--intersection $5$, 
which is the minimum self--intersection of such a linear system;
this implies that the linear system in question is exactly $|6;2^7,1^3|$. 
	
Next we consider the case $n=10$ and $r=4$, so that $h=7$;
the minimum self intersection here is $6$ 
	and the system that achieves this bound is $|4;1^{10}|$ 
	(see Theorem \ref {thm:irreduciblen+reven}, 
	and Lemma \ref{g'=0}(b)).  
	This system determines a morphism $X_{10} \longrightarrow \mathbb P^4$, 
	whose image is a surface $S'$ of degree $6$. 
	We claim that $S'$ is smooth except for $12$ non--normal simple double points 
	(transverse intersections of two smooth branches of $S'$). 
	This is an immediate consequence of the \emph{double point formula} 
	(see \cite {Cat}) and of the fact, proved in \cite [pp. 483--84]{Con}, 
	that the general internal projection of $S'$ to $\mathbb P^3$ 
	is a quintic surface with a double curve consisting of a rational normal cubic curve. 
	This in turn is the image of $X_{11}$ 
	via the linear system $|4;1^{11}|$ of dimension $3$ 
	that has also minimal self--intersection. 
	It is not difficult to directly check all this, but we will not dwell on this here.
	
	Finally, let us look at the case $n=10$ and  $r=5$.  	
	The minimum self intersection here is $8$ 
	and the system that achieves this bound is $|6;2^6,1^{4}|$ 
	(see Theorem \ref {thm:irreduciblen+rodd}, (a)). 	
	We claim that this system is very ample.
	To see this we argue as follows. 
	First consider the system $|6;2^6|$ on $X_6$. 
	This is very ample, being the double of $|3;1^6|$, 
	that is already very ample, 
	mapping $X_6$ to a smooth cubic surface in $\mathbb P^3$. 
	So $|6;2^6|$ maps $X_6$ to the $2$--Veronese image 
	of a smooth cubic surface in $\mathbb P^3$, 
	which is a surface $\Sigma$ of degree $12$ in $\mathbb P^9$ 
	with hyperplane sections of genus $4$. 
	By Green's theorem (see \cite {Gr}), 
	$\Sigma$ is $N_3$ so, in particular, it is cut out by quadrics. 
	This implies that its general internal projection 
	$\Sigma'\subset \mathbb P^8$ of degree $11$ 
	is smooth and is $N_2$, again by Green's theorem, 
	so it is cut out by quadrics. 
	So its general internal projection $\Sigma''\subset \mathbb P^7$ of degree $10$ 
	is smooth and $N_1$, again by Green's theorem, so it is still cut out by quadrics. 
	Its general internal projection $\Sigma'''\subset \mathbb P^6$ of degree $9$ is smooth. 
	Now we want to show that  the general internal projection 
	$S\subset \mathbb P^5$ is smooth, 
	which will prove our claim about the very ampleness of $|6;2^6,1^{4}|$. 
	
	Suppose this is not the case. 
	Let $p\in \Sigma'''$ be a general point. 
	Then there is a trisecant line to $\Sigma'''$ passing through $p$, 
	i.e., a line $r$ containing $p$ and cutting on $\Sigma'''$ 
	a curvilinear scheme (containing $p$) of length at least $3$. 
	Let $C\subset \mathbb P^5$ 
	be a  general hyperplane  section of $\Sigma'''$ 
	containing $r$. 
	By Bertini's theorem, $C$ is smooth. 
	The hyperplanes in $\mathbb P^5$ containing $r$ cut out on $C$, 
	which has genus $4$,  
	a linear series $g^3_d$, with $d\leq 6$. 
	This implies that $d=6$, that the $g^3_d$ is the canonical series on $C$, 
	and that $r$ is a \emph{proper trisecant}, 
	i.e., it cuts out on $\Sigma'''$ a scheme of length exactly $3$. 	
	Then $S$ would have a double point $q$ 
	(at the projection of $r$ from $p$) 
	and the projection $S'$ of $S$ from $q$ to $\mathbb P^4$ 
	would have canonical hyperplane sections. 
	This means that we would have a blow--up $X$ of $\mathbb P^2$ at $12$ points 
	(the first $10$ points are general, the remaining two are not), 
	and a linear system of the type $|6;2^6, 1^6|$ mapping $X$ to $S'$, 
	whose characteristic system is the canonical system. 
	On the other hand the canonical system is also cut out 
	on the curves of $|6;2^6, 1^6|$ by the system $|3;1^6|$ 
	and therefore this would imply that $|3;1^{12}|$ is effective. 
	But then also $|3;1^{10}|$ would be effective, a contradiction, 
	because the first blown--up 10 points are general. 													
\end{example}

We expect that linear systems on $X_n$, with $n\geq 10$, of dimension $r\geq 5$, 
with minimal self--intersection, are very ample, as it happens for $n=10$. 
Though we have explicit descriptions of these linear systems, 
this is not trivial to check in general.

\end{document}